\documentclass[parskip=half,bibliography=totoc]{scrartcl}

\usepackage[automark]{scrlayer-scrpage}								%KOMA styles					
\usepackage{amsfonts}												%Mathematics fonts
\usepackage{amsmath}												%general math tools
\usepackage{mathtools}												%General mathematics symbols
\usepackage{stmaryrd}												%Extra math symbols
\usepackage{amssymb}												%More symbols
\usepackage{extarrows}												%Extendible arrows
\usepackage{dsfont} 												%Identity matrix symbol
\usepackage{mathrsfs}												%To get mathscr
\usepackage{accents}												%Accents on math symbols
\usepackage[T1]{fontenc}											%Accents output improvement
\usepackage[utf8]{inputenc}											%Accents input improvement
\usepackage{subcaption} 											%Subfigures etc
\usepackage[top=1in,bottom=1in,left=1.25in,right=1.25in]{geometry}	%margins
\usepackage{footmisc}												%Some footnote margin thing
\usepackage{enumerate} 												%Numbered lists
\usepackage{booktabs}												%improved tables
\usepackage{authblk}												%author/affiliation formatting
\usepackage{abstract}												%for the abstract
\usepackage{tikz}

\usepackage{tikz-cd}												%Commutative Diagrams

\usepackage{amsthm}													%Theorems etc
\usepackage[colorlinks=true,citecolor=blue]{hyperref}				%hyperlinks
\usepackage[noabbrev]{cleveref}										%better cross-referencing
\crefname{lem}{lemma}{lemmata}
\crefname{prop}{proposition}{propositions}

\newcommand{\pd}[2]{\frac{\partial #1}{\partial #2}}
\renewcommand{\d}{\mathrm{d}}

\newcommand{\CH}{\mathbb{C} \mathrm{H}}
\newcommand{\N}{\mathbb{N}}
\newcommand{\Z}{\mathbb{Z}}
\newcommand{\Q}{\mathbb{Q}}
\newcommand{\R}{\mathbb{R}}
\newcommand{\C}{\mathbb{C}}
\renewcommand{\H}{\mathbb{H}}
\newcommand{\U}{\mathrm{U}}
\newcommand{\SU}{\mathrm{SU}}
\newcommand{\Sp}{\mathrm{Sp}}

\renewcommand{\Im}{\operatorname{Im}}
\renewcommand{\Re}{\operatorname{Re}}

\DeclareMathOperator{\im}{im}

\DeclareMathOperator{\Span}{span}

\DeclareMathOperator{\End}{End}
\DeclareMathOperator{\Aut}{Aut}

\newcommand{\aut}{\mathfrak{aut}}
\newcommand{\isom}{\mathfrak{isom}}
\DeclareMathOperator{\Isom}{Isom}
\DeclareMathOperator{\Lie}{Lie}

\DeclareMathOperator{\id}{id}
\newcommand{\Unit}{\mathds{1}}

\DeclareMathOperator{\Heis}{Heis}
\newcommand{\heis}{\mathfrak{heis}}

\DeclareMathOperator{\Ad}{Ad}

\DeclareMathOperator{\vol}{vol}

\newcommand{\sympquotient}{\mathbin{\hspace{-2pt}/\mkern-6mu/\hspace{-2pt}}}

\newcommand{\h}{\mathrm{H}}

\newcommand{\abs}[1]{\lvert #1 \rvert}
\newcommand{\norm}[1]{\lVert #1 \rVert}

\newcommand{\mf}[1]{\mathfrak{#1}}
\newcommand{\mc}[1]{\mathcal{#1}}

\newtheoremstyle{mythm}% name of the style to be used
{}% measure of space to leave above the theorem. E.g.: 3pt
{}% measure of space to leave below the theorem. E.g.: 3pt
{\slshape}% name of font to use in the body of the theorem
{}% measure of space to indent
{\bfseries\sffamily}% name of head font
{.}% punctuation between head and body
{ }% space after theorem head; " " = normal interword space
{}% Manually specify head
\newtheoremstyle{mydef}% name of the style to be used
{}% measure of space to leave above the theorem. E.g.: 3pt
{}% measure of space to leave below the theorem. E.g.: 3pt
{}% name of font to use in the body of the theorem
{}% measure of space to indent
{\bfseries\sffamily}% name of head font
{.}% punctuation between head and body
{ }% space after theorem head; " " = normal interword space
{}% Manually specify head

\theoremstyle{mythm}
\newtheorem{thm}{Theorem}[section]
\newtheorem{prop}[thm]{Proposition}
\newtheorem{cor}[thm]{Corollary}
\newtheorem{lem}[thm]{Lemma}
\theoremstyle{mydef}
\newtheorem{mydef}[thm]{Definition}
\newtheorem{rem}[thm]{Remark}

\newtheorem*{not*}{Notation}
\newenvironment{myproof}[1][\proofname]{
	\proof[\sffamily\upshape#1]
}{\endproof}

\newcommand{\proofclear}{\hfill \qedsymbol}

\newcommand{\hkvf}{\mathcal{X}}
\newcommand{\qkvf}{Y}
\newcommand{\hkheis}{\mathcal{U}}
\newcommand{\qkheis}{V}

\clearpairofpagestyles
\ihead[]{\headmark}
\ohead[]{\pagemark}
\deffootnote[1em]{0em}{1em}{%
	\textsuperscript{\thefootnotemark}%
}
\setfootnoterule{3em}

\newcommand\numberthis{\stepcounter{equation}\tag{\theequation}}

\newenvironment{numberedlist}{\begin{enumerate}[\upshape(i)]}{\end{enumerate}}

\title{Complete quaternionic K\"ahler manifolds with finite volume ends}
\author{V.\ Cort\'es}
\author{M.\ R\"oser}
\author{D.\ Thung}
\affil{\normalsize  Department of Mathematics \\
University of Hamburg\\
Bundesstra\ss e 55, D-20146 Hamburg, Germany}

\date{\large \today}

\begin{document}
\maketitle

\begin{abstract}
	We construct examples of complete quaternionic K\"ahler manifolds with an end of finite volume,	which are not locally homogeneous. The manifolds are aspherical with fundamental group which is up to an infinite cyclic extension a semi-direct product of a lattice in a semi-simple group with a lattice in a Heisenberg group. Their universal covering is a cohomogeneity one deformation of a symmetric space of non-compact type.

	\emph{Keywords: quaternionic K\"ahler manifolds, $c$-map, one-loop deformation, isometry groups, cohomogeneity one}\par
	\emph{MSC classification: 53C26.}
\end{abstract}

\clearpage

\tableofcontents
\clearpage

\section{Introduction}

One of the milestones of twentieth century differential geometry is Berger's classification of irreducible holonomy groups of non-locally symmetric Riemannian manifolds \cite{Ber1955}. The problem to construct non-locally symmetric Riemannian manifolds (or to show their non-existence) for each of the groups in Berger's list has been a major driving force for research in the theory of Einstein manifolds, involving a combination of Lie theory \cite{Ale1968,Ale1975}, differential geometry \cite{Bry1987,BS1989,LeB1989,LeB1991}, geometric analysis \cite{Yau1978,Joy1996,Joy1996a} and algebraic geometry \cite{Bea1983}. As a result, we have compact non-locally symmetric examples of Riemannian manifolds for all of the groups in Berger's list with exception of the group $\Sp(n)\Sp(1)$. 

Riemannian manifolds with holonomy group a subgroup of $\Sp(n)\Sp(1)$ for $n\geq 2$ are known as quaternionic K\"ahler manifolds. These manifolds are Einstein and, in fact, generalize half-conformally flat Einstein four-manifolds, which are by definition the quaternionic K\"ahler manifolds of dimension four. 

All the known compact examples of non-zero scalar curvature are locally symmetric. They fall into two classes:  symmetric quaternionic K\"ahler manifolds of compact type,  known as Wolf spaces \cite{Wol1965}, and locally symmetric spaces obtained as smooth quotients of the non-compact duals of the Wolf spaces by lattices of the isometry group \cite{Bor1963}. The manifolds in both classes are rigid \cite{LS1994,Hor1996} which shows that further examples cannot be obtained from deformation theory. Moreover,  in the case of positive scalar curvature it is even conjectured that the Wolf spaces exhaust all complete (and hence compact, by Bonnet-Myers) examples \cite{LS1994}. 

In this article we show that complete quaternionic K\"ahler manifolds of negative scalar curvature can have ends of finite volume without being locally homogeneous. Our constructions are based on the recently proven fact that the known homogeneous quaternionic K\"ahler manifolds of negative scalar curvature and higher rank can be deformed into complete quaternionic K\"ahler manifolds with an isometric cohomogeneity one action \cite[Cor.~3.19]{CST2021}. 

In this paper we focus on the above deformation $(\bar N_n,g^c)$, $c\geq 0$, for the symmetric spaces 
\begin{equation*}
	(\bar{N}_n,g^0) =\left( \frac{\SU(2,n)}{\mathrm{S}(\U(2)\times \U(n))}, g_{\mathrm{can}}\right) .
\end{equation*}
The one-parameter family $(\bar{N}_n,g^c)$, $c\geq 0$, is obtained from an indefinite hyper-K\"ahler manifold by the HK/QK correspondence. This implies that the metrics $g^c$ are quaternionic K\"ahler \cite[Cor.~1]{ACDM2015}. The completeness of $(\bar{N}_n,g^c)$ for $c\geq 0$ was shown in \cite[Cor.~15]{CDS2017}. The manifolds $(\bar{N}_n,g^c)$ are locally inhomogeneous for $c>0$ \cite[Thm.~4.8]{CST2021}.

We show in dimensions $4$ and $8$ (i.e.\ for $n=1$ and $2$) that the isometry group of  $(\bar{N}_n,g^c)$, $c\geq 0$, contains discrete subgroups which give rise to smooth quotient manifolds $X=X_\Gamma = \bar{N}/\Gamma$ with exactly two ends, see \Cref{thm:mainn=2}. The metric $g^c$ on $\bar{N}_n$ induces a complete quaternionic K\"ahler metric on $X_{\Gamma}$, which we denote by the same symbol. With respect to this metric one of the two ends is of finite volume and the other is of infinite volume, see \Cref{thm:volume}.

The fundamental groups $\Gamma$ of the eight-dimensional complete quaternionic K\"ahler ma\-nifolds $(X_\Gamma,g^c)$, $c\geq 0$, have the structure of a cyclic quotient $(\Gamma_1\ltimes \Gamma_2)/\Z$ of a semi-direct product $\Gamma_1\ltimes \Gamma_2$, where $\Gamma_1\subset \widetilde{\SU}(1,1)$ is the preimage of an arithmetic subgroup $\bar{\Gamma}_1\subset \SU(1,1)$ in the universal covering and $\Gamma_2$ is a lattice in the five-dimensional Heisenberg group. The subgroup $\Z\subset \Gamma_1\ltimes \Gamma_2$ is diagonal for $c>0$, in the sense that it has trivial intersection with each of the factors. We consider also $\bar\Gamma =\bar\Gamma_1\ltimes\Gamma_2$, which is a lattice in $\SU(1,1)\ltimes\Heis_{5}$. In the undeformed case $c=0$, we have $\Z\subset \Gamma_1$ and the above construction reduces to a semi-direct product $\Gamma =\bar\Gamma= \bar{\Gamma}_1\ltimes \Gamma_2$, of an arithmetic group $\bar{\Gamma}_1 \subset \SU(1,1)$ with a lattice $\Gamma_2 \subset \Heis_5$. In four dimensions, $\Gamma$ reduces to a lattice in the three-dimensional Heisenberg group and $X_{\Gamma}$ is diffeomorphic to $(\Heis_3/\Gamma) \times \R$.

Our results are developed in such a way that they immediately extend to higher dimensions given lattices $\bar\Gamma\subset\SU(1,n-1)\ltimes\Heis_{2n+1}$, which do always exist. In particular, the isometric action of $\SU(1,n-1)\ltimes\Heis_{2n+1}$ on a cyclic quotient of $\bar{N}_n$ is proven in all dimensions, see  \Cref{thm:SUltimesHaction}. The arithmetic part of the construction consists in finding lattices $\bar\Gamma_1$ in $\SU(1,n-1)$ which normalize a lattice in $\Heis_{2n+1}$. For $n=2$ we were able to construct co-compact lattices of this type using quaternion algebras. The lattices which we obtained can be described as follows.

Let $a,b\in \Z_{>0}$ be such that $b$ is prime and $a$ is a quadratic non-residue mod $b$, i.e.\ $a$ does not have a square root in the field $\Z_b$. Denote by $\mc O_{a,b}\subset \mf{gl}(2,\C)$ the $\Z$-span of the four matrices  
\begin{equation*}
		\Unit =
		\begin{pmatrix}
			1 & 0 \\ 0 & 1
		\end{pmatrix},\quad
		I = \sqrt a i
		\begin{pmatrix}
			0 & 1 \\ -1 & 0
		\end{pmatrix},\quad
		J=\sqrt b
		\begin{pmatrix}
			0 & 1 \\ 1 & 0
		\end{pmatrix},\quad
		K=\sqrt{ab}i 
		\begin{pmatrix}
			1 & 0 \\ 0 & -1
		\end{pmatrix}.
\end{equation*}
Then $\bar\Gamma_{1,a,b} = \mc O_{a,b}\cap\SU(1,1) = \{A\in \mc O_{a,b} : \det A = 1\}$ is a co-compact Fuchsian group, which preserves the lattice $\Gamma_2$   in $\Heis_5\cong \C^2\times\R$ generated by $\mc O_{a,b}\cdot 
\begin{psmallmatrix}
1 \\ 0
\end{psmallmatrix}
\subset \C^2$. Up to passing to a suitable finite index normal subgroup, $\bar\Gamma = \bar\Gamma_{1,a,b}\ltimes\Gamma_2$ gives a lattice of the desired  type, as long as the one-loop parameter $c$ is chosen to be a rational multiple of $\sqrt{ab}$.

The following theorem summarizes the main results about quaternionic K\"ahler manifolds of dimension $\le 8$ obtained in this paper, compare \Cref{thm:mainn=2} and \Cref{thm:volume}.
\begin{thm} 
	There exist complete locally inhomogeneous quaternionic K\"ahler manifolds of dimension $4$ and $8$ which are diffeomorphic to a product $X=\mathbb{R} \times K$, where $K$ is a compact aspherical manifold described more in detail below and the volume of one of the two ends of $X$ is finite, the other infinite. 
	\begin{itemize}
		\item For $\dim X =4$ the fiber $K$ is a nilmanifold. 
		\item For $\dim X=8$ it is a manifold with fundamental group 
		of the form $(\Gamma_1\ltimes \Gamma_2)/\mathbb{Z}$, where  $\Gamma_1$ is a lattice in $\widetilde{\SU}(1,1)$ covering a finite index normal subgroup of any of the above Fuchsian groups $\bar\Gamma_{1,a,b}\subset \SU(1,1)$, $\Gamma_2$ is a lattice in $\mathrm{Heis}_5$ and $\mathbb{Z}$ is a diagonally embedded central subgroup.  
	\end{itemize}
\end{thm}

While co-compact lattices in  $\SU(1,n-1)$ are known from complex hyperbolic geometry, see for instance \cite{Eps1987}, there is no such lattice normalizing a lattice in $\mathbb{C}^n$  for $n\ge 3$ compatible\footnote{Compatibility is defined on page \pageref{compatible:def_in_text}.}  with the Hermitian structure (and a fortiori no such lattice normalizing a compatible log-lattice in $\Heis_{2n+1}$).  Indeed, such a lattice in $\SU(1,n-1)$, $n\geq 3$, would give rise to an anisotropic Hermitian sesquilinear form over an imaginary number field. The real part thereof would be an anisotropic non-degenerate rational symmetric bilinear form in more than $4$ variables that is indefinite over $\R$, contradicting Meyer's theorem (see \cite[Cor.~2 on p.~43]{Ser73}). 

However, lattices of finite co-volume in $\SU(1,n-1)\ltimes\Heis_{2n+1}$ can be easily obtained by considering normalizers of certain lattices in $\Heis_{2n+1}$. Using suitable lattices 
\begin{equation}\label{Gamma:eq} 
	\bar\Gamma \subset \SU(1,n-1)\ltimes\Heis_{2n+1}
\end{equation} 
we construct quaternionic K\"ahler manifolds of dimension $4n$ fibering over $\R$ with locally homogeneous fibers of finite volume. The results can be summarized as follows, see \Cref{thm:finitevol} and \Cref{thm:volume}.

\begin{thm}
	In all dimensions $4n\ge 4$, there exist complete locally inhomogeneous quaternionic K\"ahler manifolds diffeomorphic to $X=\mathbb{R}\times K$, where the fibers $\{ t\} \times K$ are aspherical of finite volume with fundamental group an infinite cyclic extension of a group $\bar\Gamma$ of the form (\ref{Gamma:eq}). Infinitely many examples of such groups $\bar\Gamma$ do occur in every dimension. The domains $\{t>t_0\}\subset X$ are of finite volume for all $t_0\in \R$ while the domains $\{t<t_0\}$ are of infinite volume. 
\end{thm}

After this summary of the main results, we would like to briefly describe the structure of our paper. Our constructions are based on an explicit construction of quaternionic K\"ahler metrics of negative scalar curvature known as the one-loop deformed $c$-map. This is reviewed in \Cref{sec:prelim}. 

In \Cref{sec:action} we specialize the construction to the complete quaternionic K\"ahler $4n$-manifolds $\bar{N} = (\bar{N}_n,g^c)$, $c>0$, considered above. We prove that the simply connected group $\widetilde{\SU}(1,n-1)\ltimes \mathrm{Heis}_{2n+1}$ acts on  $\bar{N}$ by isometries and characterize its image in the isometry group. 

In \Cref{sec:quotients} we study lattices in that semi-direct product, which act with co-finite volume on the leaves of a certain codimension one foliation of $\bar N$. When the action on the leaves is co-compact, we prove that the corresponding quotients of $\bar N$ have an end of finite volume, while the other end has infinite volume.

{\bfseries Acknowledgements}

This research was partially funded by the Deutsche Forschungsgemeinschaft (DFG, German Research Foundation) under Germany's Excellence Strategy -- EXC 2121 Quantum Universe  -- 390833306. V.C. thanks Oliver Baues, Florin Belgun, and Vincent Koziarz for helpful discussions and Christoph B\"ohm for useful comments on the manuscript. We are moreover thankful to Oliver Baues and Vincent Koziarz for pointing out the references \cite{Eps1987} and \cite{Ser73}, respectively. We thank the anonymous referee for the detailed report and for making a number of useful suggestions to improve the manuscript.

\section{Preliminaries}
\label{sec:prelim}

\subsection{The (one-loop deformed) \texorpdfstring{$c$}{c}-map}

The one-loop deformed $c$-map is a string theory construction  \cite{RSV2006} which, to a given projective special K\"ahler manifold $\bar M$, associates a one-parameter family of quaternionic K\"ahler manifolds. A geometric proof of the quaternionic K\"ahler property was given in \cite[Theorem 5]{ACDM2015} based on an extension \cite{ACM2013} of Haydys' HK/QK correspondence \cite{Hay2008}. These constructions were elegantly recovered in \cite{MS2015} based on Swann's twist approach \cite{Swa2010}, as we will now describe. For more information see~\cite{CST2021,ACDM2015,MS2015} and references therein.

We begin by recalling \cite{Fre1999,ACD2002} the notion of a projective special K\"ahler manifold, which is defined in terms of a corresponding conical affine special K\"ahler manifold.

\begin{mydef}
	A pseudo-K\"ahler manifold $(M,g,J,\omega)$ endowed with a flat, torsion-free connection $\nabla$ which satisfies $\nabla\omega=0$ and $\d^\nabla J=0$, where $J$ is viewed as a $TM$-valued one-form, is called an affine special K\"ahler (ASK) manifold. 
	
	An ASK manifold $(M,g,J,\nabla)$ endowed with a vector field $\xi$ is called conical, or a CASK manifold, if 
	\begin{numberedlist}
		\item $\{\xi,J\xi\}$ generate a principal (that is, free and proper) $\C^*$-action,
		\item $g|_{\langle \xi,J\xi\rangle}$ is negative-definite while $g|_{\langle \xi,J\xi\rangle^\perp}$ is positive-definite, and
		\item $\nabla \xi=\nabla^\mathrm{LC}\xi=\id_{TM}$, where $\nabla^\mathrm{LC}$ is the Levi-Civita connection of $g$.
	\end{numberedlist}
	The vector field $\xi$ is called the Euler field.
\end{mydef}

The $\C^*$-action on a CASK manifold casts it as the total space of a $\C^*$-principal bundle. The base of this bundle is what we call a projective special K\"ahler manifold.

\begin{mydef}
	Let $(M,g,J,\nabla,\xi)$ be a CASK manifold. Then $\bar M\coloneqq M/\C^*$ is called a projective special K\"ahler (PSK) manifold.
\end{mydef}

That every PSK manifold is indeed K\"ahler follows from the fact that we can construct $\bar M$ as a K\"ahler quotient. Given a CASK manifold $(M,g,J,\nabla,\xi)$, $-J\xi$ generates an $S^1$-action which preserves the (pseudo-)K\"ahler structure. Its Hamiltonian function is $\mu=-\frac{1}{2}g(\xi,\xi)$, whose level sets intersect the orbits of the $\R_+$-action generated by the Euler field exactly once. Choosing the level $\frac{1}{2}$, we then find $M\sympquotient S^1=\mu^{-1}\big(\frac{1}{2}\big)/S^1\cong M/\C^*=\bar M$. Since $g|_{\langle \xi,J\xi\rangle^\perp}$ is positive-definite, so is the induced K\"ahler metric on $\bar M$.

By definition, $\bar M$ arises from a CASK manifold $(M,g_M,J_M,\omega_M,\nabla,\xi)$. The ASK structure on $M$ gives rise to a pseudo-hyper-K\"ahler structure on its cotangent bundle.

\begin{thm}[\cite{CFG1989,Fre1999,ACD2002}]
	Let $(M,g_M,J_M,\omega_M,\nabla)$ be an ASK manifold. Then $N\coloneqq T^*M$ carries a pseudo-hyper-K\"ahler structure. With respect to the splitting $TN\cong \pi^*TM\oplus \pi^*T^*M$ induced by $\nabla$, the hyper-K\"ahler metric $g$ and complex structures $I_k$, $k=1,2,3$, are given by the expressions
	\begin{equation*}
		g=
		\begin{pmatrix}
			g_M & 0 \\ 0 & g^{-1}_M
		\end{pmatrix}\hspace{1.25cm}
		I_1=
		\begin{pmatrix}
			J_M & 0 \\ 0 & J^*_M
		\end{pmatrix}\hspace{1.25cm}
		I_2=
		\begin{pmatrix}
			0 & -\omega_M^{-1}\\
			\omega_M & 0
		\end{pmatrix}\hspace{1.25cm}
		I_3=I_1I_2,
	\end{equation*}
	where $g^{-1}_M$ and $J_M^*$ denote the natural induced structures on $T^*M$, $\omega_M$ is regarded as an isomorphism $TM\to T^*M$, and we have omitted pullbacks throughout to simplify notation. \proofclear
\end{thm}

This pseudo-hyper-K\"ahler structure, which is known as the rigid $c$-map structure, is a natural generalization of the standard pseudo-hyper-K\"ahler structure on $\H^{n,1}\cong\C^{n,1}\oplus (\C^{n,1})^*$. 

In the conical case, we additionally have a natural circle action, generated by the vector field $-J_M\xi$.

\begin{mydef}
	We call a vector field $Z\in \mf X(N)$ on a (pseudo-)hyper-K\"ahler manifold $(N,g,\omega_k)$, $k=1,2,3$, a rotating Killing field if it satisfies $L_Zg=0$, $L_Z\omega_1=0$, $L_Z\omega_2=\omega_3$ and $L_Z\omega_3=-\omega_2$.
\end{mydef}

Using coordinates adapted to the CASK structure~\cite{ACD2002}, a straightforward computation shows:

\begin{prop}[\cite{ACM2013}]\label{prop:ACM}
	Let $(M,g_M,J_M,\nabla,\xi)$ be a CASK manifold, and endow $N=T^*M$ with the rigid $c$-map structure. Then $Z=-\widetilde{J_M\xi}\in \mf{X}(N)$, where the tilde denotes the $\nabla$-horizontal lift, is a rotating Killing field. Moreover, it is $\omega_1$-Hamiltonian, with Hamiltonian function $f_Z=-\frac{1}{2}g(Z,Z)$.\proofclear
\end{prop}

In summary, the CASK structure on $M$ induces a pseudo-hyper-K\"ahler structure equip\-ped with a free circle action, generated by an $\omega_1$-Hamiltonian rotating Killing field, on its cotangent bundle. This is precisely the input required for the so-called HK/QK correspondence~\cite{Hay2008,ACM2013,ACDM2015}, which produces a one-parameter family of quaternionic K\"ahler manifolds $(\bar N,g^c)$, $c\geq 0$, endowed with a nowhere-vanishing Killing field, out of it. The fact that, for the data described in \Cref{prop:ACM}, the resulting metric is precisely the one-loop deformed $c$-map metric of \cite{RSV2006} was proven in \cite{ACDM2015}.

The twist construction \cite{Swa2010} provides a duality between manifolds endowed with a distinguished, closed two-form and a vector field which is Hamiltonian with respect to it. In the case at hand, the relevant two-form is $\omega_\h\coloneqq \omega_1+\d \iota_Z g$, and the vector field is the rotating Killing field $Z$ introduced above. For any $c\in \R$, the function $f_\h=f_Z+g(Z,Z)=\frac{1}{2}g(Z,Z)-\frac{1}{2}c$ is a Hamiltonian for $Z$ with respect to $\omega_\h$, i.e.~satisfies $\d f_\h=-\iota_Z\omega_\h$. The triple $(Z,\omega_\h,f_\h)$ makes up the so-called twist data for the HK/QK correspondence. In the following, we will assume that $\omega_\h$ is an integral two-form, i.e.~has integer periods.

Exploiting the integrality, we pass to a principal $S^1$-bundle $\pi_N:P\to N$ endowed with a principal $S^1$-connection $\eta$ with curvature $\omega_\h$. The next step is to lift $Z\in \mf X(N)$ to a vector field on $P$ which preserves $\eta$; such a lift is provided by $Z_P\coloneqq \tilde Z+\pi^*_N f_\h X_P$, where $\tilde Z\in \mf X(P)$ is the $\eta$-horizontal lift of $Z$ and $X_P\in \mf X(P)$ is generator of the principal circle action on $P$. It is known that, after imposing an integrality condition on the additive constant $c\in \R$ which appears in the definition of $f_\h$, one may assume that $Z_P$ generates a circle action~\cite{Swa2010}. If the circle action on $N$ is free, as in the case at hand, we may moreover assume that the lifted action is free as well. The quotient of $P$ by this lifted circle action is then a smooth manifold; this quotient realizes $P$ as the total space of a second circle bundle $\pi_{\bar N}:P\to \bar N$. Since the two principal circle actions on $P$ commute, $\bar N$ inherits a twisted circle action, generated by $(\pi_{\bar N})_*X_P$.

The power of the twist construction lies in the fact that it enables one to use geometric structures on $N$ to induce twisted structures on $\bar N$, assuming that they are invariant under the initial circle action. This is done by exploiting the double fibration structure (see~\eqref{eq:twistfibrations}), which induces pointwise identifications of the tangent spaces to both $N$ and $\bar N$ with the $\eta$-horizontal subspaces of the tangent spaces to $P$.

\begin{equation}\label{eq:twistfibrations}
	\begin{aligned}
		\begin{tikzpicture}
			\draw (1,1.5) node{$N$};
			\draw[<->,dashed] (1.9,1.5) -- (4,1.5)
				node[pos=0.5, anchor=north]{twist};
			\draw (3,3.5) node{$P$};
			\draw[->] (2.75,3.25) -- (1.25,1.85) node[pos=0.5,anchor = south east]{$\pi_N$};
			\draw[->] (3.25,3.25) -- (4.75,1.85) node[pos=0.5,anchor = south west]{$\pi_{\bar N}$};
			\draw (5,1.5) node{$\bar N$};
		\end{tikzpicture}
	\end{aligned}
\end{equation}

In particular, the quaternionic structure bundle of the hyper-K\"ahler manifold $N$, which is preserved by the rotating Killing field $Z$, endows $\bar N$ with an almost quaternionic structure. However, the hyper-K\"ahler metric $g$ does not twist to a quaternionic K\"ahler metric on $\bar N$. In order to obtain a quaternionic K\"ahler structure, we must first perform a so-called elementary deformation~\cite{MS2015}. 

A (pseudo-)Riemannian metric $h$ on $N$ is called an elementary deformation if it takes the following form: 
\begin{equation*}
	h=a g+ b \sum_{\mu=0}^3\alpha_\mu\otimes \alpha_\mu
\end{equation*}
where, after setting, $\omega_0\coloneqq g$, we use the notation $\alpha_\mu=\iota_Z\omega_\mu$, and $a,b\in C^\infty(N)$ are nowhere-vanishing functions. Assuming that $Z$ is nowhere-vanishing, the second term is proportional to the restriction of $g$ to the quaternionic span $\H Z=\langle Z,I_1 Z,I_2 Z,I_3 Z\rangle$ of $Z$. Thus, we can think of an elementary deformation as arising from the composition of a conformal scaling with an independent conformal scaling along $\H Z$.

Building on earlier work of several authors, Macia and Swann showed that there is an essentially unique combination of elementary deformation and twist data which yield a (pseudo-)quaternionic K\"ahler metric:

\begin{thm}[\cite{Hay2008,ACM2013,ACDM2015,MS2015}]\label{thm:twistHKQK}
	Let $(N,g,\omega_k,Z,f_Z)$ be a (pseudo-)hyper-K\"ahler manifold endowed with an $\omega_1$-Hamiltonian, rotating Killing field. Assume that its Hamiltonian function $f_Z$ is nowhere-vanishing and let $k_1,k_2\in \R\setminus\{0\}$. Then, with respect to the twist data $(Z,\omega_\h,f_\h)$, where
	\begin{equation*}
		\omega_\h=k_1(\omega_1+\d\iota_Z g) \qquad \qquad \qquad f_\h=k_1(f_Z+g(Z,Z))
	\end{equation*}
	the elementary deformation
	\begin{equation*}
		g_\h\coloneqq k_2\bigg(\frac{1}{f_Z}g+\frac{1}{f_Z^2}\sum_\mu \alpha_\mu\otimes \alpha_\mu\bigg)
	\end{equation*}
	twists to a (pseudo-)quaternionic K\"ahler metric on $\bar N$. These are the only choices of elementary deformation and twist data which lead to a (pseudo-)quaternionic K\"ahler structure. Moreover, when the initial (pseudo-)hyper-K\"ahler manifold arises from a CASK manifold in the manner described above, the resulting quaternionic K\"ahler metrics are positive-definite, and of negative scalar curvature. \proofclear
\end{thm}

The constants $k_1,k_2$ do not influence the local geometry of the resulting quaternionic K\"ahler manifold, and will be set to $1$ in the following. However, the additive constant $c$ parametrizing the freedom in the choice of the Hamiltonian function $f_Z$ does play an important role as can be seen from the explicit appearance of $f_Z$ in the formula for $g_\h$. 

In our setup, where $N$ is the cotangent bundle of a CASK manifold, we have $f_Z=-\frac{1}{2}(g(Z,Z)+c)$. Denoting the resulting metrics by $g^c$, we have in summary obtained a one-parameter family $(\bar N,g^c)$, $c\geq 0$, of quaternionic K\"ahler manifolds out of a PSK manifold $\bar M$. The metric $g^0$ is known as the undeformed $c$-map metric (or Ferrara--Sabharwal metric); the metrics $g^c$, $c>0$, are called (one-loop) deformed $c$-map metrics.

We will rely on the following completeness results:

\begin{thm}[\cite{CHM2012,CDS2017}]\label{thm:completeness}
	Let $(\bar N,g^c)$, $c\geq 0$, be the image of the PSK manifold $\bar M$ under the deformed $c$-map. Then $(\bar N,g^0)$ is complete if and only if $\bar M$ is complete. If, additionally, $\bar M$ arises from the (supergravity) $r$-map, or has regular boundary behavior, then $(\bar N,g^c)$ is complete for all $c\geq 0$. \proofclear
\end{thm}

Another fact to keep in mind is:

\begin{prop}[\cite{CDS2017}]\label{prop:differentcisometric}
	Let $(\bar N,g^c)$, $c\geq 0$, be the image of the PSK manifold $\bar M$ under the deformed $c$-map. Then, for arbitrary $c_1,c_2>0$, $(\bar N,g_{c_1})$ is locally isometric to $(\bar N,g_{c_2})$. \proofclear
\end{prop}

The undeformed $c$-map metric, however, is generally distinct, as we will soon demonstrate in explicit examples.
  
\subsection{Symmetries and the \texorpdfstring{$c$}{c}-map}
\label{sec:symmetries}

In this section we review what is known about the symmetries of the quaternionic K\"ahler metrics constructed via the $c$-map and its one-loop deformation. In particular, we outline how one may obtain Killing vector fields from infinitesimal automorphisms of the initial PSK manifold.

It is well-known that the $4n$-dimensional quaternionic K\"ahler manifold $(\bar N,g^c)$ that arises from the $c$-map applied to a PSK manifold $\bar M$ of complex dimension $n-1$ admits an isometric action of a $2n+1$-dimensional Heisenberg group. In the examples considered in this paper, where $\bar M$ is a so-called PSK domain, we can think of $(\bar N,g^c)$ as a bundle over $\bar M$ such that the Heisenberg group preserves fibers and acts on them freely. In particular, the orbits are submanifolds of the fibers of codimension one. For $c=0$, this can be improved: A one-dimensional solvable extension of the Heisenberg group acts isometrically, freely and transitively on the fibers, casting $(\bar N,g^0)$ as a bundle of Lie groups. A detailed description, including an explicit formula for the action with respect to a standard coordinate system on the fibers, is given in~\cite[Lem.~2.20]{CST2021}.

Given this large group of fiber-preserving isometries, it is natural to ask whether (deformed) $c$-map metrics admit additional isometries that cover non-trivial diffeomorphisms of the initial PSK manifold $\bar M$. In particular, if the PSK structure of $\bar M$ admits non-trivial automorphisms, we may ask if these can be lifted to symmetries of its image under the $c$-map. 

In order to provide a precise answer, we must formalize what an automorphism of a PSK manifold is. Keeping in mind the extrinsic definition of this class of manifolds, the following is natural:

\begin{mydef}
	Let $(M,g,J,\nabla,\xi)$ be a CASK manifold, and $\bar M$ the corresponding PSK manifold. An automorphism of the CASK structure on $M$ is a diffeomorphism of $M$ which preserves the pseudo-K\"ahler structure and the flat connection $\nabla$, and commutes with the $\C^*$-action induced by $\{\xi,J\xi\}$. An automorphism of the PSK manifold $\bar M$ is a diffeomorphism of $\bar M$ which is induced by a CASK automorphism of the corresponding CASK manifold $M$.
\end{mydef}

The matter of lifting PSK automorphisms to the image $(\bar N,g^c)$ under the $c$-map was taken up in~\cite{CST2021}. By definition, every one-parameter group of PSK automorphisms lifts to a one-parameter group of CASK automorphisms of the corresponding CASK manifold $M$. These CASK automorphisms can in turn be lifted canonically to $N=T^*M$, using the pullback on one-forms. It turns out that the canonical lift of a CASK automorphism to $N$ preserves all data on $N$ that we have discussed, i.e.~the pseudo-hyper-K\"ahler structure and twist data $(Z,\omega_\h,f_\h)$. Moreover, every one-parameter family of canonically lifted CASK automorphisms is $\omega_\h$-Hamiltonian, and there is a canonical choice of Hamiltonian function. In fact, there exists an equivariant moment map for the lifted action of $\Aut M$:

\begin{prop}
	Let $(M,g,\omega,\nabla,\xi)$ be an CASK manifold. Then an equivariant moment map for the (canonically lifted) action of $\Aut M$ on $N$, with respect to $\omega_\h$, is given by $\mu:N\to \mf g^*$ where $\langle\mu,X\rangle\eqqcolon\mu_X=\frac{1}{2}\big(g(Z,X)+\nu_X\big)$. Here, $\nu_X(\alpha)=(\omega^{-1})_{\pi(\alpha)}(\alpha\circ\nabla V^X,\alpha)$, $\alpha\in T^*M$, with $V^X$ the fundamental vector field corresponding to $X$.
\end{prop}
\begin{myproof}
	The expression given here for the Hamiltonian function corresponding to $X\in \aut M$ agrees (though the notation differs) with the results of~\cite{CST2021}, so we will not verify (again) that it indeed provides a Hamiltonian.  
	One then checks the equivariance, i.e.~that $\mu_X(g\cdot \alpha)=\mu_{\Ad_{g^{-1}}X}(\alpha)$ for all $g\in \Aut M$ and $\alpha\in T^*M$, by working through the definitions of $\mu_X$ and $\nu_X$. 
\end{myproof}

\begin{thm}[\cite{CST2021}]\label{thm:liftingsymmetries}
	Let $\aut(\bar M)$ denote the algebra of infinitesimal automorphisms of the PSK manifold $\bar M$, and $\isom(\bar N,g^c)$ the algebra of Killing fields of its image $(\bar N,g^c)$ under the (deformed) $c$-map. Then there exists an injective, linear map $\aut(\bar M)\to \isom(\bar N,g^c)$ for every $c\geq 0$.
\end{thm}
\begin{myproof}[Sketch of Proof]
	 After lifting an infinitesimal PSK automorphism to the corresponding CASK manifold, and then to its cotangent bundle $N$, let us denote the generator of the lifted one-parameter family by $X$, and its $\omega_\h$-Hamiltonian by $f_X$. Then we certainly have $L_X g_\h=0$, where $g_\h$ is the elementary deformation from~\Cref{thm:twistHKQK}. However, since the twist of a Killing field is generally not a Killing field with respect to the twisted metric, we must modify $X$ to obtain a Killing field of the quaternionic K\"ahler metric $g^c$. This aim is achieved by $X-\frac{f_X}{f_\h}Z$, whose twist $Y$ indeed satisfies $L_Y g^c=0$. Note that $Y$ implicitly depends on the deformation parameter $c$ through the appearance of $f_\h$.
\end{myproof}
\begin{rem}\label{rem:lifts}\leavevmode
	\begin{numberedlist}
		\item The proof shows that the entire procedure can be applied just as well to an infinitesimal CASK automorphism of $M$ that does not arise as a lift of an infinitesimal PSK automorphism of $\bar M$. Such infinitesimal automorphisms can only exist in case the CASK manifold is flat (cf.~\cite{CST2021}), in which case the generator of the distinguished circle action is an infinitesimal automorphism. The series of examples discussed in this paper is of this type.
		\item The principal circle action on $P$ pushes down to a distinguished, isometric circle action on $(\bar N,g^c)$ which commutes with the Killing fields constructed by means of \Cref{thm:liftingsymmetries}.		
	\end{numberedlist}
\end{rem}

The map $\aut(\bar M)\to\isom(\bar N,g^c)$ of \Cref{thm:liftingsymmetries} is in general not a Lie algebra homomorphism. Upon twisting, the commutator of vector fields (generically) picks up an additional term which is proportional to the vector field generating the distinguished circle action on $\bar N$, yielding a one-dimensional central extension: 

\begin{thm}[\cite{CST2021}]\label{thm:liftedalgebraicstructure}
	The Lie algebra $\aut(\bar M)$ induces a Lie algebra of Killing fields on $(\bar N,g^c)$ which is isomorphic to a one-dimensional central extension of $\aut(\bar M)$.\proofclear
\end{thm}

\begin{rem}
If the central extension is trivial, then we may arrange matters such that we do obtain a subalgebra of $\isom(\bar N,g^c)$ isomorphic to $\aut(\bar M)$. This is the case in the series of examples on which we focus in this work.
\end{rem}

The Killing fields constructed via \Cref{thm:liftingsymmetries} are linearly independent of the generators of the action of the $2n+1$-dimensional Heisenberg group. Putting them together and taking into account \Cref{rem:lifts}, we have constructed an algebra of Killing fields of $(\bar N,g^c)$ whose dimension is $\dim\aut(M)+2n+1$, where $\dim_\C M=n$. For $c=0$, this is improved to $\dim \aut(M)+2n+2$. Under the assumptions of \Cref{thm:completeness}, we know that $(\bar N,g^c)$ is complete for all $c\geq 0$ and can integrate these Killing fields to obtain a group of isometries of the (deformed) $c$-map metric.

\section{The \texorpdfstring{$c$}{c}-map applied to complex hyperbolic spaces}
\label{sec:action}
\subsection{The infinitesimal action}

We will now focus on an important series of examples. Consider $\C^n$, equipped with the flat pseudo-K\"ahler structure induced by the indefinite Hermitian form $h=-\d z^0\otimes \d \bar z^0+\sum_{a=1}^{n-1}\d z^a\otimes\d \bar z^a$. Since the Levi-Civita connection is flat, we may regard it as endowing $\C^n$ with the structure of an ASK manifold. The standard $\C^*$-action by scalar multiplication endows the $\C^*$-invariant domain $M_n=\{z\in \C^n\mid -\abs{z^0}^2 + \sum_{a=1}^{n-1}\abs{z^a}^2<0 \}$ with the structure of a CASK manifold. The K\"ahler quotient by $\U(1)\subset \C^*$ is nothing but complex hyperbolic space $\CH^{n-1}$ equipped with the (symmetric) Bergman metric
\begin{equation}\label{eq:Bergmanmetric}
	g_{\CH^{n-1}}=\frac{1}{1-\norm{X}^2}\Bigg(\sum_{a=1}^{n-1} \abs{\d X^a}^2+\frac{1}{1-\norm{X}^2}\bigg|\sum_{a=1}^{n-1}\bar X^a \d X^a\bigg|^2\Bigg)
\end{equation}
The indefinite unitary group $\U(1,n-1)$ acts on $M_n$ by automorphisms of the CASK structure and projects down to yield a transitive action on $\CH^{n-1}$, which is therefore a homogeneous PSK manifold.

It follows directly from the general theory developed in~\Cref{sec:symmetries} that the image $(\bar N_n,g^0)$ of $\CH^{n-1}$ under the undeformed $c$-map is a homogeneous quaternionic K\"ahler manifold (of negative scalar curvature), and that the deformed metrics $(\bar N_n,g^c)$, $c>0$, are complete and of cohomogeneity at most one. An explicit expression for the metrics $g^c$, with respect to a global coordinate system $(X^a, w^k,\tilde\phi,\rho )\in \C^{n-1}\times \C^n\times \R\times \R_{>0}$, where $a=1,\dots,n-1$ and $\norm{X}^2=\sum_a \abs{X^a}^2<1$, and $k=0,\dots,n-1$, was given in \cite{CDS2017}:
\begin{align*}
	g^c&=\frac{\rho+c}{\rho}g_{\CH^{n-1}}+\frac{1}{4\rho^2}\frac{\rho+2c}{\rho+c}\d \rho^2\\
	&\quad +\frac{1}{4\rho^2}\frac{\rho+c}{\rho+2c}
	\bigg(\d \tilde\phi 
	-4\Im\bigg(\bar w^0 \d w^0-\sum_{a=1}^{n-1}\bar w^a \d w^a\bigg)
	+\frac{2c}{1-\norm{X}^2}\Im \bigg(\sum_{a=1}^{n-1} \bar X^a \d X^a\bigg)\bigg)^2\\\numberthis\label{eq:fullmetric}
	&\quad -\frac{2}{\rho}\bigg(\d w^0 \d \bar w^0-\sum_{a=1}^{n-1}\d w^a \d \bar w^a\bigg)
	+\frac{\rho+c}{\rho^2}\frac{4}{1-\norm{X}^2}	\Big|\d w^0+\sum_{a=1}^{n-1}X^a\d w^a\Big|^2
\end{align*}

Note that the case $n=1$, where we apply the $c$-map to a single point, is exceptional. With the convention $\sum_{a=1}^0=0$, the expression \eqref{eq:fullmetric} remains valid and greatly simplifies. The resulting quaternionic K\"ahler four-manifold is known as the universal hypermultiplet in the physics literature. Despite the qualitative difference with the case $n\geq 2$, our results in this paper remain valid if $n=1$, unless explicitly indicated.

As a smooth manifold, $\bar N_n$ is just a copy of $\R^{4n}$. However, its Riemannian structure is interesting. In fact, it is known that the undeformed $c$-map metric casts $\bar N_n$ as the non-compact symmetric space $\frac{\SU(n,2)}{\mathrm{S}(\U(n)\times \U(2))}$. Moreover, it was proven in~\cite{CST2021} that $(\bar N_n,g^c)$ is of cohomogeneity one. In summary, the family $(\bar N_n,g^c)$, $c\geq 0$ is a one-parameter deformation of a symmetric space through cohomogeneity one complete quaternionic K\"ahler metrics. 

To develop a more concrete understanding of these manifolds, we will now show in detail how the general theory works out in these examples. Our main goal is to derive explicit formulas for the Killing fields obtained by means of \Cref{thm:liftingsymmetries}. The starting point is the (standard) action of $\U(1,n-1)$ on $M_n$ by CASK automorphisms, or rather the infinitesimal action of its Lie algebra $\mf u(1,n-1)$. Passing to $N_n=T^*M_n=M_n\times \C^n\subset \C^n\times \C^n$, the canonically lifted vector field $\hkvf_A$ corresponding to $A\in \mf u(1,n-1)$ is given by $\hkvf_A(z,w)=(Az+\overline{Az},-A^\top w-\overline{A^\top w})$, where we identify the tangent spaces with the underlying vector space.

The twisting two-form $\omega_\h$ is 
\begin{equation*}
	\omega_\h=\frac{i}{2}\bigg(\d z^0\wedge\d \bar z^0-\d w^0\wedge \d \bar w^0 
	-\sum_{a=1}^{n-1}\big(\d z^a\wedge\d\bar z^a-\d w^a\wedge\d \bar w^a\big)\bigg),
\end{equation*}
so that the Hamiltonian corresponding to $\hkvf_A$ is given by $f_{\hkvf_A}(z,w)=-\frac{1}{2}\Im\big(\langle Az,z\rangle+\langle Aw,w\rangle\big)$, where $\langle z_1,z_2\rangle=-z_1^0\bar z_2^0+\sum_{a=1}^{n-1}z_1^a\bar z_2^a$. The rotating Killing field $Z$, which is the horizontal lift of the (negative) generator of the standard $\U(1)$-action on the first factor, is given by $Z(z,w)=(-iz,0)$. Its Hamiltonian is 
\begin{equation*}
	f_\h=f_Z+g(Z,Z)=\frac{1}{2}(g(Z,Z)-c)=\frac{1}{2}(\langle z,z\rangle -c).
\end{equation*} 
 (Recall that for the particular vector field $Z$ we denoted by $f_Z=-\frac12 (g(Z,Z)+c)$ its Hamiltonian with respect to $\omega_1$. That is why 
we are using the notation $f_H$ for the $\omega_H$ Hamiltonian.)
Following the proof of \Cref{thm:liftingsymmetries}, we now know that the twist of $\hkvf_A^\h\coloneqq \hkvf_A-\frac{f_{\hkvf_A}}{f_\h}Z$ will be a Killing field on $(\bar N_n,g^c)$. 

In order to derive explicit expressions for these twisted vector fields, we first construct a circle bundle $P$, equipped with a connection $\eta$ with curvature $\omega_\h$. Since $\omega_\h$ is exact, $P$ is trivial and we may choose $\eta=\d s+\frac{1}{2}\iota_\Xi\omega_\h$, 
where $s$ is the periodic coordinate parametrizing $S^1=\R/2\pi \Z$ and $\Xi$ is the position vector field
\begin{equation*}
	\Xi=\sum_{k=0}^{n-1}\bigg(z^k \pd{}{z^k} + \bar z^k \pd{}{\bar z^k} +w^k \pd{}{w^k}+\bar w^k \pd{}{\bar w^k}\bigg).
\end{equation*}
Explicitly, 
\begin{equation*}
	\eta=\d s+\frac{i}{4}\Big(z^0\d \bar z^0-\bar z^0\d z^0-w^0\d \bar w^0+\bar w^0\d w^0
	-\sum_{a=1}^{n-1}\big(z^a\d \bar z^a-\bar z^a\d z^a-w^a \d \bar w^a+\bar w^a\d w^a\big)\Big).
\end{equation*}
We lift $Z$ to the vector field $Z_P$, which generates the circle action defining the principal circle bundle $P\to \bar N_n$. Recall that $Z_P=\tilde Z+f_\h \pd{}{s}$, where the tilde denotes the $\eta$-horizontal lift and we regard $f_\h$ as a function on $P=N_n\times S^1$. Regarding $Z$ as a vector field on $P$, a short computation shows that 
\begin{equation*}
	Z_P=Z-\frac{c}{2}\pd{}{s}=-i\sum_{k=0}^{n-1}\bigg(z^k \pd{}{z^k}-\bar z^k \pd{}{\bar z^k}\bigg)-\frac{c}{2}\pd{}{s}.
\end{equation*}
Similarly, for $A\in\mf{u}(1,n-1)$, one may check that $\eta(\hkvf_A^\h)=-\frac{f_{\hkvf_A}}{f_\h}\frac{c}{2}$, so that
\begin{equation*}
	\tilde \hkvf_A^\h=\hkvf_A^\h+\frac{f_{\hkvf_A}}{f_\h}\frac{c}{2}\pd{}{s}
	=\hkvf_A-\frac{f_{\hkvf_A}}{f_\h}Z_P.
\end{equation*}
We now want to push this vector field down to $\bar N_n$. To this end, we realize $\bar N_n$ as a submanifold of $P$, using a global slice for the action obtained by integrating $Z_P$, which is given by $(z,w,s)\mapsto (e^{-it}z,w,s-\frac{c}{2}t)$, $t\in \R$. The flow is not periodic for arbitrary $c\geq 0$, but for $c\in 2\Z$ we have a circle action covering the circle action on $N_n$ (this is the previously mentioned quantization condition). Since the local geometry of the metrics $g^c$, $c>0$, is independent of $c$ (cf.~\Cref{prop:differentcisometric}) we will only consider the values $c=0$ and $c=2$. The submanifold $\bar N_n\coloneqq\{(z,w,s)\in N_n\times S^1\mid \arg z^0=0\}$ then provides a global slice for the action, and can therefore be identified as the twist manifold. 

The final ingredient required is an identification of the horizontal tangent spaces to $P$ with the tangent spaces of $\bar N_n$, given by an $S^1$-principal connection with respect to the circle action generated by $Z_P$, namely $\theta=\frac{1}{f_\h}\eta$. Now we have to write $\tilde \hkvf_A^\h$ in the form $\tilde\hkvf_A^\h=\qkvf_A-\theta(\qkvf_A)Z_P$, where $\qkvf_A$ is tangent to $\bar N_n\subset P$ (at every point in $\bar N_n$). $\qkvf_A$ is then be the twist of $\hkvf_A^\h$, whose $\theta$-horizontal lift indeed agrees with $\eta$. Since $\theta(Z_P)=1$ and $\theta(\tilde \hkvf_A^\h)=0$ by construction, we already have
\begin{equation*}
	\tilde \hkvf_A^\h=\hkvf_A-\theta(\hkvf_A)Z_P,
\end{equation*}
but $\hkvf_A$ will not be tangent to $\bar N_n$ in general. To remedy this, we have to consider a vector field of the form $\hkvf_A+\varphi Z_P$, where $\varphi$ is a function. This vector field restricts to $\bar N_n$ if and only if it satisfies $\d(\arg z^0)(\hkvf_A+\varphi Z_P)=\d(\arg z^0)(\hkvf_A)-\varphi=0$, where we used that $\d(\arg z^0)(Z_P)=-1$, showing that we must have $\varphi=\d(\arg z^0)(\hkvf_A)$. The twisted vector field $\qkvf_A$ on $\bar N_n$ is therefore given by 
\begin{equation}\label{eq:twist}
	\qkvf_A\coloneqq \hkvf_A+\d(\arg z^0)(\hkvf_A)Z_P.
\end{equation}

Now we turn to the algebra $\mf u(1,n-1)$ of infinitesimal CASK automorphisms of $M_n$, which we consider as a real form of $\mf{gl}(n,\C)$. Writing $A\in \mf{gl}(n,\C)$ in block form with respect to the splitting $\C^n = \C e_0 \oplus e_0^\perp$ (where $\{e_k\}_{k=0,\dots,n-1}$ is an orthonormal basis with respect to the indefinite Hermitian form), i.e. 
\begin{equation}\label{eq:decomp}
	A = 
	\begin{pmatrix}
		\lambda & v^\top  \\ w & B
	\end{pmatrix},
\end{equation}
with $\lambda\in \C,v,w\in \C^{n-1}$ and $B\in \mf{gl}(n-1,\C)$, we will think of $\mf u(1,n-1)$ as the fixed-point set of the anti-linear involutive automorphism 
\begin{equation*}
	\sigma: \mf{gl}(n,\C)\to\mf{gl}(n,\C),\qquad A^\sigma\coloneqq\sigma(A)\coloneqq -I\bar{A}^\top I,
\end{equation*}
where $I =
\begin{psmallmatrix}
	-1 & 0 \\ 0 & \Unit_{n-1}
\end{psmallmatrix}$. 
We may now write any $A\in\mf{gl}(n,\C)$ uniquely in the form 
\begin{equation*}
	A = \Re_{\sigma}(A) + i \Im_{\sigma}(A), 
\end{equation*}
where 
\begin{equation*}
	\Re_{\sigma}(A) = \frac{1}{2}(A + A^\sigma), \qquad \Im_{\sigma}(A) = \frac{1}{2i}(A-A^\sigma)
\end{equation*}
are both elements of $\mf{u}(1,n-1)$. Note that $A^\sigma=\Re_\sigma(A)-i\Im_\sigma(A)$.

The Lie algebra $\mf{sl}(n,\C)$ is generated by the matrices 
\begin{equation}\label{eq:UaUasigma}
	U_a = 
	\begin{pmatrix}
		0 & e_a^\top  \\ 0 & 0 
	\end{pmatrix},\qquad 
	U_a^\sigma\coloneqq\sigma(U_a) = 
	\begin{pmatrix}
		0 & 0 \\ e_a & 0 
	\end{pmatrix}, \qquad a = 1,\dots,n-1.
\end{equation}
We observe that 
\begin{equation}\label{eq:LieUaUbsigma}
	[U_a,U_b] = 0, \qquad [U_a^\sigma,U_b^\sigma]= 0, \qquad 
	[U_a,U_b^\sigma] = 
	\begin{pmatrix}
		\delta_{ab} & 0 \\ 0 & -e_be_a^\top 
	\end{pmatrix}.
\end{equation}
Thus, 
\begin{equation*}
	\{U_a,U_a^\sigma, [U_a,U_b^\sigma]: a,b = 1,\dots,n-1\}
\end{equation*}
is a basis of $\mf{sl}(n,\C)$. From this, together with the central element $C = i\Unit_n$, they form a basis for $\mf{gl}(n,\C)$, which is therefore generated (as a complex Lie algebra) by $\{C,U_a,U_a^\sigma\colon a=1,\dots,n-1\}$.

The infinitesimal action of $\mf u(1,n-1)$ on $M_n\subset \C^n$ by CASK automorphisms extends complex-linearly to an infinitesimal action $\mf{gl}(n,\C)\to \mf X^\C(M_n)$ by complex vector fields. Since both canonically lifting and the twisting procedure extend to complex vector fields without problems, we can complexify the above discussion of $\mf u(1,n-1)$ to obtain a linear map $\alpha^\C:\mf{gl}(n,\C)\to \mf X^\C(\bar N_n)$.

\begin{prop}\label{prop:alphagln}
	The linear map $\alpha^\C: \mf{gl}(n,\C)\to \mf{isom}(\bar N_n,g^c)^\C$, $A\mapsto \qkvf_A$ defines an infinitesimal action of $\mf{gl}(n,\C)$ by complex Killing fields on $(\bar N_n,g^c)$. It maps $\mf u(1,n-1)$ into $\isom(\bar N_n,g^c)$ and is given on generators by  
	\begin{equation}\label{eq:alphagln}
		\begin{aligned}
			\qkvf_C&\coloneqq \alpha^\C(C)=-i\sum_{k=0}^{n-1}\bigg(w^k\pd{}{w^k}-\bar w^k \pd{}{\bar w^k}\bigg)-2c\pd{}{\tilde\phi},\\
			\qkvf_a&\coloneqq\alpha^\C(U_a)= \pd{}{\bar X^a}-X^a\sum_{b=1}^{n-1}X^b \pd{}{X^b}
			-w^0\pd{}{w^a}-\bar w^a\pd{}{\bar w^0} + ic X^a \pd{}{\tilde\phi},\\
			\bar\qkvf_a&=\alpha^\C(U_a^\sigma),
		\end{aligned}
	\end{equation}
	where $(X^a\coloneqq \frac{z^a}{z^0}, w^k, \tilde\phi\coloneqq 4s, \rho\coloneqq 2f_Z=\abs{z^0}^2-\sum_{a=1}^{n-1}\abs{z^a}^2-c)$, $a=1,\dots, n-1$, $k=0,\dots ,n-1$, are global coordinates on 	$\bar N_n\cong B^{n-1}\times \R^{2n+1}\times \R_{>0}$.  Here $B^{n-1}\subset \C^{n-1}$ denotes the unit ball.
\end{prop}

\begin{cor}\label{Cor: KillingYaVa}
	The real vector fields $Y_C,\Re(Y_a), \Im(Y_a)\in \mf X(\bar N_n)$ generate a subalgebra of $\isom(\bar N_n,g^c)$ isomorphic to $\mf u(1,n-1)$. \proofclear
\end{cor}

\begin{myproof}[Proof of \Cref{prop:alphagln}]	
	First observe that it follows immediately from the fact that $\alpha^\C$ satisfies $\alpha^\C(A^\sigma)=\overline{\alpha^\C(A)}$ for every $A\in \mf{gl}(n,\C)$ that $\alpha^\C(\mf u(1,n-1))\subset \mf{isom}(\bar N_n,g^c)$. 
	
	Our next aim is to prove \eqref{eq:alphagln}, for which we use \eqref{eq:twist}. Denoting the canonical lift of the fundamental vector field (on $M_n$) corresponding to $U_a$ by $\hkvf_a$, we obtain
	\begin{equation*}
		\hkvf_a=z^a\pd{}{z^0}+\bar z^0 \pd{}{\bar z^a}-w^0 \pd{}{w^a}-\bar w^a\pd{}{\bar w^0}.
	\end{equation*}
	Since $\d(\arg z^0)=-\frac{i}{2}\big(\frac{\d z^0}{z^0}-\frac{\d \bar z^0}{\bar z^0}\big)$, $\d(\arg z^0)(\hkvf_a)=-\frac{i}{2}\frac{z^a}{z^0}$ and the twist of $\hkvf_a$ is
	\begin{align*}
		\qkvf_a&=\hkvf_a-\frac{i}{2}\frac{z^a}{z^0}Z_P\\
		&=z^a\pd{}{z^0}+\bar z^0\pd{}{\bar z^a}-w^0 \pd{}{w^a}-\bar w^a\pd{}{\bar w^0}
		-\frac{1}{2}\frac{z^a}{z^0}\sum_{k=0}^{n-1}\bigg(z^k \pd{}{z^k}-\bar z^k \pd{}{\bar z^k}\bigg)+\frac{ic}{4}\frac{z^a}{z^0}\pd{}{s}.
	\end{align*}
	Expressing this vector field in terms of the conventional (global) coordinate functions $(X^a\coloneqq \frac{z^a}{z^0}, w^k, \tilde\phi\coloneqq 4s, \rho\coloneqq 2f_Z=\abs{z^0}^2-\sum_{a=1}^{n-1}\abs{z^a}^2-c)$ on $\bar N_n$, we obtain $Y_a$ as in \eqref{eq:alphagln}. Applying the same procedure to the central element $C=i\Unit_n$, we obtain $Y_C$ as in \eqref{eq:alphagln}. Finally, the third line of \eqref{eq:alphagln} follows from $\alpha^\C(U_a^\sigma) =\overline{\alpha^\C(U_a)}$.
	
	Next, we explicitly compute the commutators of the complex Killing fields $\qkvf_C$, $\qkvf_a$, $\bar\qkvf_a$ and compare them with the structure constants of $\mf{gl}(n,\C)$ with respect to the basis $\{C,U_a,U_a^\sigma,[U_a,U_b^\sigma]\colon a,b=1,\dots,n-1\}$.
	
	Since $C=i\Unit_n$ lies in the center of $\mf{gl}(n,\C)$, we start by verifying that $\qkvf_C$ lies in the center of $\Im\alpha^\C\subset \mf X^\C(\bar N_n)$. Since $C\in \mf u(1,n-1)$, $\qkvf_C$ is real and 
	\begin{align*}
		[\qkvf_C,Y_a] &= -i\Bigg[\sum_{k=0}^{n-1}\bigg(w^k\pd{}{w^k}- \bar w^k\pd{}{w^k}\bigg)+2c\pd{}{\tilde\phi}\, , \, 
		\pd{}{\bar X^a} - \sum_j X^aX^j\pd{}{X^j}+ ic X^a \pd{}{\tilde\phi} \Bigg]\\
		&\quad\ + i\Bigg[\sum_{k=0}^{n-1}\bigg(w^k\pd{}{w^k}- \bar w^k\pd{}{w^k}\bigg)\, ,\,  
		w^0\pd{}{w^a} + \bar w^a\pd{}{w^0} \Bigg]\\
		&= i\Bigg[\sum_{k=0}^{n-1}\bigg(w^k\pd{}{w^k}- \bar w^k\pd{}{w^k}\bigg) ,  
		w^0\pd{}{w^a} + \bar w^a\pd{}{w^0} \Bigg]= 0. 
	\end{align*}
	It follows that $\qkvf_C$ commutes with $\bar Y_a$ and thus also with the commutators $[Y_a,\bar Y_b]$ for all $a,b$. This means that $\qkvf_C$ commutes with every vector field in the image of $\alpha^\C$, as desired. 
	
	Next, since $[U_a,U_b] =0$  we have to check that $[Y_a,Y_b] = 0$ for all $a,b$:
	\begin{align*}
		[Y_a,Y_b] &= \Bigg[\pd{}{\bar X^a} - \sum_j X^aX^j\pd{}{X^j}  + ic X^a \pd{}{\tilde\phi}\, , \, 
		\pd{}{\bar X^b} - \sum_k X^bX^k\pd{}{X^k}  + ic X^b \pd{}{\tilde\phi}\Bigg] \\ 
		& \quad \ + \bigg[w^0\pd{}{w^a}+ \bar w^a\pd{}{\bar w^0}\, , \,  
		w^0\pd{}{w^b}+ \bar w^b\pd{}{\bar w^0}\bigg]\\
		&= \Bigg[- \sum_j X^aX^j\pd{}{X^j}  + ic X^a \pd{}{\tilde\phi}\, , \,   
		- \sum_k X^bX^k\pd{}{X^k}  + ic X^b \pd{}{\tilde\phi}\Bigg] \\ 
		&=   \sum_{j,k} \Big[X^aX^j\pd{}{X^j}\, ,\, 
		X^bX^k\pd{}{X^k}\Big]  
		- ic\bigg(\sum_j X^aX^j\delta_{jb} - \sum_k X^bX^k\delta_{ka}\bigg) \pd{}{\tilde\phi}\\
		&= 0.
	\end{align*}
	This relation immediately implies $[\bar Y_a,\bar Y_b]=0$. 
	
	Finally, we show that $\alpha^\C([U_a,U_b^\sigma])=-[Y_a,\bar Y_b]$; note the sign, which arises due to the fact that the map sending a Lie algebra element to its fundamental vector field is a Lie algebra anti-homomorphism in our setting. We start by computing the commutator $[Y_a,\bar Y_b]$:
	\begin{align*}
		[Y_a,\bar Y_b] &= \Bigg[\pd{}{\bar X^a} - \sum_j X^aX^j\pd{}{X^j}  + ic X^a \pd{}{\tilde\phi}\, , \, \pd{}{X^b} - \sum_j \bar X^b\bar X^j\pd{}{\bar X^j}  - ic \bar X^b \pd{}{\tilde\phi}\Bigg] \\ 
		&\quad \  + \bigg[w^0\pd{}{w^a}+ \bar w^a\pd{}{\bar w^0}\, , \,  \bar w^0\pd{}{\bar w^b}+  w^b\pd{}{ w^0}\bigg]\\
		&= -\sum_j\bigg((\delta_{ab}\bar X^j + \delta_{aj}\bar X^b)\pd{}{\bar X^j} - (\delta_{ab}X^j + \delta_{bj}X^a)\pd{}{X^j}\bigg) -2ic\delta_{ab}\pd{}{\tilde\phi} \\ 
		& \quad \ +  w^0\delta_{ab}\pd{}{w^0} - w^b\pd{}{w^a} + \bar w^a \pd{}{\bar w^b} - \delta_{ab}\bar w^0\pd{}{\bar w^0}\\
		&= \delta_{ab}\Bigg(\sum_j\bigg(X^j\pd{}{X^j}-\bar X^j\pd{}{\bar X^j} \bigg)  +  w^0\pd{}{w^0}- \bar w^0\pd{}{\bar w^0} -2ic\pd{}{\tilde\phi}\Bigg)\\
		& \quad \ + X^a\pd{}{X^b} - \bar X^b\pd{}{\bar X^a} + \bar w^a \pd{}{\bar w^b}  - w^b\pd{}{w^a}.
	\end{align*}
	On the other hand, the canonical lift of the fundamental vector field corresponding to $[U_a,U_b^\sigma]$ is
	\begin{equation*}
		\hkvf_{[U_a,U_b^\sigma]}= \delta_{ab}\bigg(z_0\pd{}{z_0}-\bar{z}_0\pd{}{\bar z_0} -w_0\pd{}{w_0}+\bar{w}_0\pd{}{\bar w_0}\bigg) - z_a\pd{}{z_b}+\bar{z}_b\pd{}{\bar z_a} + w_b\pd{}{w_a}-\bar{w}_a\pd{}{\bar w_b}.
	\end{equation*}
	Upon twisting, we indeed obtain
	\begin{align*}
		\alpha^\C([U_a,U_b^\sigma])&=-\delta_{ab}\Bigg(\sum_j\bigg(X^j\pd{}{X^j}-\bar X^j\pd{}{\bar X^j} \bigg)  +  w^0\pd{}{w^0}- \bar w^0\pd{}{\bar w^0} -2ic\pd{}{\tilde\phi}\Bigg)\\
		& \quad \ - X^a\pd{}{X^b} + \bar X^b\pd{}{\bar X^a} - \bar w^a \pd{}{\bar w^b} + w^b\pd{}{w^a}.
	\end{align*}
	This completes the proof.
\end{myproof}

The \Cref{prop:alphagln} provides an (infinitesimal) action of the unitary Lie algebra $\mf{u}(1,n-1)$ on $\bar N_n$. As explained in \Cref{sec:symmetries}, there is an additional action of a $(2n+1)$-dimensional Heisenberg group. As we will now show, this action also has a natural interpretation in terms of the twist construction of $\bar N_n$. 

On the hyper-K\"ahler side of the twist correspondence, translations of the $w$-coordinates define an action of $\R^{2n}$ on $N_n$. The action preserves the pseudo-hyper-K\"ahler structure and twist data $(Z,\omega_\h,f_\h)$, and is $\omega_\h$-Hamiltonian. Thus, we may apply the procedure outlined in the proof of \Cref{thm:liftingsymmetries} to turn its infinitesimal generators into Killing fields on $\bar N_n$. 

We will once again work with the complexified generators $\pd{}{w^k}$, $k=0,1\dots,n-1$, and their conjugates. A natural choice of $\omega_\h$-Hamiltonian function for $\pd{}{w^k}$ is $f_{w^k}=\pm \frac{i}{2}\bar w^k$, where the sign is positive for $k=0$ and negative otherwise. Carrying out the exact same procedure as outlined for the vector fields arising from the underlying CASK manifold, we now obtain Killing fields on $\bar N_n$ by twisting the modified vector fields $\hkheis_k=\pd{}{w^k}-\frac{f_{w^k}}{f_\h}Z$. Horizontally lifting, we find
\begin{align*}
	\tilde\hkheis_k&=\hkheis_k-\eta(\hkheis_k)\pd{}{s}
	=\hkheis_k-\frac{1}{2}\big(\Xi(f_{w^k})-\frac{f_{w^k}}{f_\h}\Xi(f_\h) \big)\pd{}{s}\\
	&=\hkheis_k+\frac{f_{w^k}}{2f_\h}(f_\h+c)\pd{}{s}=\pd{}{w^k}+\frac{f_{w^k}}{2}\pd{}{s}-\frac{f_{w^k}}{f_\h}Z_P
\end{align*}
where we used that $\Xi(f_{w^k})=f_{w^k}$ and $\Xi(f_\h)=2f_\h+c$ in passing to the second line. We note that the first two terms determine a vector field $\qkheis_k$ on $P$ which is tangent to $\bar N_n$ (at every point in $\bar N_n$) and satisfies $\tilde \hkheis_k=\qkheis_k-\theta(\qkheis_k)Z_P$, which means that $\qkheis_k$ restricts to the sought-after Killing field on $\bar N_n$. After substituting $\tilde\phi=4s$, we have:
\begin{equation}\label{eq:DefVa}
	\qkheis_a = \pd{}{w^a} -  i\bar w^a\pd{}{\tilde\phi}, \qquad \qkheis_0 =  \pd{}{w^0} +i\bar w^0 \pd{}{\tilde\phi}. 
\end{equation}
where $a=1,\dots,n-1$. The complex conjugate vector fields $\bar V_k$, $k=0,\dots,n-1$, naturally arise in the same way from $\pd{}{\bar w^k}$.

Now let $e_0,f_0,e_a,f_a$, $a=1,\dots, n-1$, be the standard basis of $\R^{2n}$. We define the one-dimensional central extension $\heis_{2n+1}$ of $\R^{2n}$ by setting
\begin{equation}\label{eq:commutators}
	[e_k,e_l]=0, \quad [f_k,f_l]=0,\quad 
	[e_k,f_l]=\bigg(\delta_{k0}\delta_{l0}-\sum_{a=1}^{n-1}\delta_{ka}\delta_{la}\bigg)T
\end{equation} 
for every $k,l=0,1,\dots,n-1$, where $T$ denotes the generator of the center. Complexifying and extending the Lie bracket complex-bilinearly, we obtain $\heis_{2n+1}^\C$. 

\begin{prop}\label{prop:heisaction}
	The vector fields $\qkheis_k$ and $\bar\qkheis_k$, $k=0,\dots,n-1$, together with $\pd{}{\tilde\phi}$, generate an infinitesimal action of a complexified Heisenberg algebra $\heis_{2n+1}^\C$ by complex Killing fields on $(\bar N_n,g^c)$ which sends $\heis_{2n+1}$ into $\isom(\bar N_n,g^c)$.
\end{prop}
\begin{myproof}	
	We define a (complex) infinitesimal action, i.e.~a Lie algebra anti-homomorphism $\beta^\C:\heis_{2n+1}^\C\to \mf X^\C(\bar N)$, by setting $\beta^\C(E_k)=V_k$ and $\beta^\C(\bar E_k)=\bar V_k$, where $k=0,1,\dots,n-1$ and $E_k\coloneqq e_k-if_k$, as well as $\beta(T)=\pd{}{\tilde\phi}$. The fact that $\beta^\C$ defines an infinitesimal action follows from 
	\begin{equation*}
		[V_k,\bar V_l] = -2i\bigg(\delta_{k0}\delta_{l0}-\sum_{a=1}^{n-1}\delta_{ka}\delta_{la}\bigg)\pd{}{\tilde\phi}
	\end{equation*}
	and the vanishing of all other Lie brackets, which is easily deduced from \eqref{eq:DefVa}. The vector fields $V_k$, $\bar V_k$ are complex Killing fields by construction. Their real and imaginary parts, which generate $\R^{2n}\subset \heis_{2n+1}$, are real Killing fields on $\bar N_n$. The general theory of the twist construction also dictates that $\pd{}{\tilde\phi}$, which is the restriction of the generator of the principal circle action of the bundle $P\to N_n$, is Killing. Indeed, the quaternionic K\"ahler metric $g^c$ arises by pushing down the pullback of the elementary deformation $g_\h$, which is certainly invariant under the principal circle action.
\end{myproof}

The next logical step is to study how the generators of the infinitesimal (complex) $\mf{gl}(n,\C)$- and $\heis_{2n+1}^\C$-actions interact with one another. To determine this interaction we recall that, on the hyper-K\"ahler side of the HK/QK correspondence, the vector spaces spanned by the $w$-coordinates are the fibers of $N_n=T^*M_n$. Therefore, $\mf u(1,n-1)$ naturally acts on them via the dual of the defining representation. 

Passing to the quaternionic K\"ahler side, these same $w$-coordinates span $\R^{2n}\subset \heis_{2n+1}$. We therefore expect that $\mf u(1,n-1)$ acts on this subspace via the dual representation. We extend this to an action on all of $\heis_{2n+1}$ by stipulating that $\mf u(1,n-1)$ acts trivially on the center. This defines a representation of $\mf u(1,n-1)$
by derivations of $\heis_{2n+1}$ and therefore yields a semi-direct product structure $\mf u(1,n-1)\ltimes \heis_{2n+1}$. Complexifying and extending the Lie bracket complex-bilinearly, we obtain a complex Lie algebra $\mf{gl}(n,\C)\ltimes \heis_{2n+1}^\C$.

Recall that $\{U_a,U_a^\sigma,[U_a,U_b^\sigma]\colon a,b=1,\dots,n-1\}$, together with $C$, form a basis of $\mf{gl}(n,\C)$. We complete to a basis of $\mf{gl}(n,\C)\ltimes\heis_{2n+1}^\C$ by adding the basis elements $T$, $E_k$, $\bar E_k$, $k=0,1,\dots,n-1$ of $\heis_{2n+1}^\C$ (as in the proof of \Cref{prop:heisaction}). We now extend the infinitesimal action $\alpha^\C:\mf{gl}(n,\C)\to \isom(\bar N_n,g^c)^\C$ to an infinitesimal action $\alpha^\C:\mf{gl}(n,\C)\ltimes \heis_{2n+1}^\C\to \isom(\bar N_n,g^c)^\C$ in the manner prescribed in the proof of \Cref{prop:heisaction}, i.e.~by setting
\begin{equation}\label{eq:Defalphaltimes}
	\alpha^\C(E_k)=V_k,\quad\alpha^\C(\bar E_k)=\bar V_k,\quad \alpha^\C(T)=\pd{}{\tilde\phi}.
\end{equation}
where $E_k=e_k-if_k$.

\begin{prop}
	The linear map $\alpha^\C\colon \mf{gl}(n,\C)\ltimes\heis_{2n+1}^\C\to\isom(\bar N_n,g^c)^\C$ as defined by \Cref{prop:alphagln} and \Cref{eq:Defalphaltimes} is an infinitesimal action on $(\bar N,g^c)$ by complex Killing fields which sends $\mf{u}(1,n-1)\ltimes\heis_{2n+1}$ into $\isom(\bar N_n,g^c)$.
\end{prop}
\begin{myproof}
	It is clear from the preceding discussion that $\alpha^\C$ indeed maps into $\isom(\bar N_n,g^c)^\C$, and that the image of $\mf{u}(1,n-1)\ltimes\heis_{2n+1}$ consists of real Killing fields of $(\bar N_n,g^c)$. All that remains is to verify that $\alpha^\C$ defines a Lie algebra anti-homomorphism. 
	
	By \Cref{prop:alphagln,prop:heisaction}, we need only check those Lie brackets that involve one element of $\mf{gl}(n,\C)$ and one element of $\heis_{2n+1}^\C$. In $\mf{gl}(n,\C)\ltimes \heis^\C_{2n+1}$, the following brackets can be calculated from \eqref{eq:UaUasigma} and the fact that, by definition of the semi-direct product structure, the bracket of $\mf{gl}(n,\C)\ltimes\heis_{2n+1}^\C$ evaluated on $A\in\mf{u}(1,n-1)$ and $v\in\R^{2n}$ is just $[A,v] = -A^\top v$ (where $A^\top $ is identified with a real $2n\times2n$-matrix). This prescription yields
	\begin{align*}
		[C,E_k] = -iE_k,&\qquad  [C,\bar E_k] = i\bar E_k,\\
		[U_a,E_k]=-\delta_{k0}E_a, & \qquad [U^\sigma_a,E_k]=-\delta_{ka}E_0,\\
		[U_a,\bar E_k]=-\delta_{ka}\bar E_0, & \qquad 
		[U^\sigma_a,\bar E_k]=-\delta_{k0}\bar E_a
	\end{align*}
	The remaining brackets involving one element of $\mf{gl}(n,\C)$ and one of $\heis_{2n+1}^\C$ can be deduced from these via the Jacobi identity. We now compare to the brackets of the corresponding vector fields on the quaternionic K\"ahler manifold, using the formulas \eqref{eq:alphagln} and \eqref{eq:DefVa}. They yield the expected results, namely
	\begin{align*}
		[Y_C,V_k] = iV_k,&\qquad [Y_C,\bar V_k]= -i\bar V_k,\\ 
		[Y_a, V_k] = \delta_{k0}V_a,&\qquad [\bar Y_a,V_k] = \delta_{ka}V_0,\\
		[Y_a,\bar V_k] = \delta_{ka}\bar V_0,&\qquad [\bar Y_a,\bar V_k] = \delta_{k0}\bar V_a.
	\end{align*}
	This confirms that $\alpha^\C$ defines an infinitesimal action.
\end{myproof}

This shows that, despite the $c$-dependence of the Killing fields constructed via the HK/QK correspondence, the isomorphism type of the Lie algebra of Killing fields is independent of $c$. The $c$-dependence arises only once we consider the ($c$-dependent) infinitesimal action $\alpha:\mf u(1,n-1)\ltimes \heis_{2n+1}\to \isom(\bar N_n,g^c)$ obtained by restricting $\alpha^\C$.

\subsection{\texorpdfstring{$\mf{u}(1,n-1)\ltimes\heis_{2n+1}$}{The image of the infinitesimal action} as an algebra of Killing fields}

The goal of this section is to provide a characterization of $\mf g\coloneqq\Im\alpha$ as a subalgebra of the Killing algebra of $\isom(\bar N_n,g^c)$, see \Cref{prop:char}. We start by determining the stabilizer of a distinguished point:

\begin{lem}\label{lem:stabilizer}
	Consider for fixed $\rho_0>0$ the point $n_0 =(0,0,0,\rho_0)\in \bar N_n$. Then the stabilizer of $n_0$ in $\mf g$ is given by 
	\begin{equation*}
		\mf g_{n_0} = \Span\bigg\{Y_C+2c\pd{}{\tilde\phi}, \Re([Y_a,\bar Y_b]), \Im([Y_a,\bar Y_b]) + 2c\delta_{ab}\pd{}{\tilde\phi} \colon a,b=1,\dots,n-1\bigg\}.
	\end{equation*}
	In particular, $\mf g_{n_0}$ is isomorphic to $\mf{u}(1)\times\mf{u}(n-1)$. 
\end{lem}
\begin{myproof}
	Recall from the proof of \Cref{prop:alphagln} the formula 
	\begin{equation}\label{eq:commexpression}
		\begin{aligned}
			[Y_a,\bar Y_b] 	&= \delta_{ab}\Bigg(\sum_{d=1}^{n-1}\bigg(X^d\pd{}{X^d}-\bar X^d\pd{}{\bar X^d} \bigg)  +  w^0\pd{}{w^0}- \bar w^0\pd{}{\bar w^0} -2ic\pd{}{\tilde\phi}\Bigg)\\
			& \quad \ + X^a\pd{}{X^b} - \bar X^b\pd{}{\bar X^a} + \bar w^a \pd{}{\bar w^b}  - w^b\pd{}{w^a}.
		\end{aligned}
	\end{equation}
	From this we read off 
	\begin{equation*}
		[Y_a,\bar Y_b](n_0)= -2ic\delta_{ab} \pd{}{\tilde\phi}.
	\end{equation*}
	Similarly, $Y_C(n_0)=-2c\pd{}{\tilde\phi}$; the inclusion of the right-hand side in the left-hand side now follows. The other inclusion can be seen directly from \eqref{eq:alphagln}, \eqref{eq:DefVa} and \eqref{eq:Defalphaltimes}. 
	
	To see that this subalgebra is isomorphic to $\mf{u}(1)\times\mf{u}(n-1)$, consider its preimage under the (injective!) infinitesimal action $\alpha$, which is
	\begin{equation}\label{eq:abstractstabilizer}
		\mf g'=\Span\big\{C+2cT,\Re([U_a,U_b^\sigma]),\Im([U_a,U_b^\sigma])+2c\delta_{ab}T
		\colon a,b=1,\dots,n-1\big\}
	\end{equation}
	Since $\heis_{2n+1}\subset \mf u(1,n-1)\ltimes \heis_{2n+1}$ is an ideal, the natural projection map $\pi_1:\mf u(1,n-1)\ltimes \heis_{2n+1}\to \mf u(1,n-1)$ is a homomorphism of Lie algebras. Its restriction to $\mf g'$ is injective since $\alpha(\heis_{2n+1})\cap \mf g_{n_0}=\{0\}$, and from \eqref{eq:abstractstabilizer} one reads off that 
	\begin{equation*}
		\pi_1(\mf g')=\Span\big\{C,\Re([U_a,U_b^\sigma]),\Im([U_a,U_b^\sigma])\colon a,b=1,\dots,n-1\big\},
	\end{equation*}
	which is nothing but the standard embedding $\mf u(1)\times \mf u(n-1)\subset \mf u(1,n-1)$.
\end{myproof}

Using this, we can determine the subalgebra consisting of Killing fields tangent to the submanifold $H\coloneqq \{X=0,\rho=\rho_0\}\subset \bar N_n$, which is acted on simply transitively and isometrically by $\Heis_{2n+1}$. Since the isometry group of the symmetric space $(\bar N_n,g^0)$ is well-understood, we focus on $(\bar N_n,g^c)$, $c>0$.

\begin{lem}\label{tgH:lem}
	Let $c>0$ and let $K$ be a Killing field on $(\bar N_n, g^c)$ which is tangent to $H$ at every point contained in $H$. Then there exists an element of the Lie subalgebra $\alpha((\mf{u}(1)\times\mf{u}(n-1))\ltimes\heis_{2n+1})\subset\mf{g}$ 
	which coincides with $K$ along $H$. 
\end{lem}
\begin{myproof}
	Clearly, $K|_H$ is a Killing field on $(H,g_H^c)$ where we wrote $g_H^c$ for the metric on $H$ induced by $g^c$. From \eqref{eq:fullmetric}, we find
	\begin{equation*}
		\begin{aligned} 
			g^c_H &=\frac{1}{2\rho_0}\sum_{a=1}^{n-1}\big|\d w^a\big|^2 + \bigg(\frac{\rho_0+c}{\rho_0^2}-\frac{1}{2\rho_0}\bigg)\big|\d w^0\big|^2\\
			& \quad+ \frac{1}{4\rho_0^2}\frac{\rho_0+c}{\rho_0+2c}\bigg(\d\tilde\phi 
			-2i\bigg(\bar w^0\d w^0 - w^0\d\bar w^0 
			-\sum_{a=1}^{n-1}(\bar w^a\d w^a-w^a\d\bar w^a)\bigg)\bigg)^2,
		\end{aligned}
	\end{equation*}
	which is a left-invariant metric on the Heisenberg group $\Heis_{2n+1}\cong H$ by construction. A result of Wilson~\cite{Wil1982} asserts that the Killing algebra of $(H,g_H^c)$ can be constructed as a semi-direct product $\mf{s}\ltimes\heis_{2n+1}$, where $\mf{s}$ is the Lie algebra of skew-symmetric derivations of $\heis_{2n+1} = T_{(0,0,0)}H$.
	
	Note that 
	\begin{equation*}
		g^c_H(0,0,0) = \frac{1}{2\rho_0}\sum_{a=1}^{n-1}\abs{\d w^a}^2
		+ \bigg(\frac{\rho_0+c}{\rho_0^2}-\frac{1}{2\rho_0}\bigg)\abs{\d w^0}^2
		+ \frac{1}{4\rho_0^2}\frac{\rho_0+c}{\rho_0+2c}\d\tilde\phi^2
	\end{equation*} 
	is diagonal, and recall, see \eqref{eq:commutators}, that the bracket on $\heis_{2n+1}$ is given by $[v_1,v_2]=\omega(v_1,v_2)\pd{}{\tilde\phi}$, where  $\omega=\frac{i}{2}\big(\d w^0\wedge \d\bar w^0 - \sum_{a=1}^{n-1}\d w^a\wedge \d\bar w^a\big)$.	
	
	Now let $A\in \End(\heis_{2n+1})$ be a skew-symmetric derivation, so that in particular 
	\begin{equation}\label{omega_eq}
		\omega(v_1,v_2)A\bigg(\pd{}{\tilde\phi}\bigg) 
		=\big([Av_1,v_2]+[v_1,Av_2]\big)\pd{}{\tilde\phi},\quad \mbox{for all}\quad v_1,v_2\in  \R^{2n}.
	\end{equation}
	For $v\in\R^{2n}$ we write $Av = \tilde Av + l_A(v)\pd{}{\tilde\phi}$ and $A\big(\pd{}{\tilde\phi}\big) = \lambda_A\pd{}{\tilde\phi}$ (every derivation preserves the center) with $l_A\in (\R^{2n})^*, \lambda_A\in \R$. Since $g^c_H$ is diagonal at $(0,0,0)$, skew-symmetry implies that $\lambda_A = 0$, $l_A= 0$. Then 
	equation (\ref{omega_eq}) simplifies to
	\begin{equation*}
		0=\omega(\tilde Av_1,v_2) + \omega(v_1,\tilde Av_2),
	\end{equation*}
	that is $\tilde{A}\in\mf{sp}(2n,\R)$.
	
	Introducing real coordinates $\{\tilde\zeta_k,\zeta^k\}$ by writing $w^0 =\frac{1}{2}(\tilde \zeta_0+i\zeta^0)$ and $w^a=\frac{1}{2}(\tilde \zeta_a-i\zeta^a)$, we have $\omega = \frac{1}{2} \sum_{k=0}^{n-1} \d\zeta^k\wedge \d\tilde\zeta_k$ and
	\begin{equation*}
		g^c_H = \frac{1}{8\rho_0}\bigg(\sum_{a=1}^{n-1}(\d\zeta^a)^2 + (\d\tilde\zeta^a)^2\bigg) +
		\frac{1}{4}\bigg(\frac{\rho_0+c}{\rho_0^2}-\frac{1}{2\rho_0}\bigg)(\d\zeta^0)^2 + (\d\tilde\zeta^0)^2)   +
		\frac{1}{4\rho_0^2}\frac{\rho_0+c}{\rho_0+2c}\d\tilde\phi^2.
	\end{equation*}
	Write $\tilde A$ in block form with respect to the splitting $\R^{2n} = \R^2 \oplus \R^{2n-2}$.  
	\begin{equation*}
		\tilde A = \begin{pmatrix}
			A_1 & A_2 \\ A_3 & A_4
		\end{pmatrix}
	\end{equation*}
	The Gram matrices of $g^c_H$ (restricted to $\R^{2n}=\Span\{\tilde\zeta_k,\zeta^k\}$) and $\omega$ are 
	\begin{equation*}
		G = \frac{1}{4}
		\begin{pmatrix} \big(\frac{\rho_0+c}{\rho_0^2}-\frac{1}{2\rho_0}\big)\Unit_2 & 0\\
			0 & \frac{1}{2\rho_0}\Unit_{2n-2}  
		\end{pmatrix}, \qquad 
		\Omega = 
		\frac{1}{2}\begin{pmatrix}
			J_2 & 0 \\ 0 & J_{2n-2}
		\end{pmatrix},
	\end{equation*}
	where $J_2= %
	\begin{psmallmatrix}
		0 & -1 \\ 1 & 0
	\end{psmallmatrix}$ %
	and $J_{2n-2}$ consists of $(n-1)$ copies of $J_2$ on its diagonal. Skew-symmetry of $A$ implies $\tilde A^\top G + G\tilde A = 0$ and $\tilde A\in \mf{sp}(2n,\R)$ implies that $\tilde A^\top \Omega + \Omega\tilde A = 0$. This means that 
	\begin{align*}
		\begin{pmatrix}
			A_1^\top  & A_3^\top  \\ A_2^\top  & A_4^\top 
		\end{pmatrix}\!
		\begin{pmatrix} 
			\big(\frac{\rho_0+c}{\rho_0^2}-\frac{1}{2\rho_0}\big)\Unit_2 & 0\\
			0 & \frac{1}{2\rho}\Unit_{2n-2}  
		\end{pmatrix}\! \! +\! \!
		\begin{pmatrix} 
			\big(\frac{\rho_0+c}{\rho_0^2}-\frac{1}{2\rho_0}\big)\Unit_2 & 0\\
			0 & \frac{1}{2\rho_0}\Unit_{2n-2}  
		\end{pmatrix}\!
		\begin{pmatrix}
			A_1 & A_2 \\ A_3 & A_4
		\end{pmatrix} \!&= 0,\\
		\begin{pmatrix}
			A_1^\top  & A_3^\top  \\ A_2^\top  & A_4^\top 
		\end{pmatrix}
		\begin{pmatrix}
			J_2 & 0 \\ 0 & J_{2n-2}
		\end{pmatrix} + 
		\begin{pmatrix}
			J_2 & 0 \\ 0 & J_{2n-2}
		\end{pmatrix}
		\begin{pmatrix}
			A_1 & A_2 \\ A_3 & A_4
		\end{pmatrix}\!&= 0.
	\end{align*}
	This implies 
	\begin{align*}
		A_1+A_1^\top  &= 0\\
		A_1^\top J_2+J_2 A_1 &= 0\\
		A_4+A_4^\top  &= 0\\
		A_4^\top J_2 + J_2 A_4 &=0 \\
		\bigg(\frac{\rho_0+c}{\rho_0^2}-\frac{1}{2\rho_0}\bigg)A_2^\top  + \frac{1}{2\rho_0} A_3 &= 0\\
		A_2^\top J_2 + J_{2n-2}A_3 &= 0
	\end{align*}
	The last two lines imply $A_2 = 0$ and $A_3=0$, since the parameter $c>0$. Thus, $\tilde A$ is of the form 
	\begin{equation*}
		\tilde A = \begin{pmatrix}
			A_1 & 0 \\ 0 & A_4
		\end{pmatrix}
	\end{equation*}
	with $A_1$ and $A_4$ skew-symmetric and complex linear. Put differently, we have shown that $A\in\mf{u}(1)\times\mf{u}(n-1)$, and thereby that $\mf s\cong \mf u(1)\times \mf u(n-1)$. 
	
	The stabilizer (in $\mf g=\Im \alpha$) of the point $(0,0,0,\rho_0)$, determined in \Cref{lem:stabilizer}, consists of Killing fields that are tangent to $H$ and span a subalgebra isomorphic to $\mf u(1)\times \mf u(n-1)$. Combining this with $\alpha(\heis_{2n+1})$, we see that $\alpha((\mf u(1)\times \mf u(n-1))\ltimes \heis_{2n+1})\cong \mf s\ltimes \heis_{2n+1}$, so the full algebra of Killing fields of $H$ arises from $\alpha$ as claimed.
\end{myproof}

The image of $\alpha^\C$, which we denote by $\mf{g}^\C$, contains the subalgebra 
\begin{equation*}
	\mf{h} = \Span_\C\Big\{Y_a, V_a, V_0, \pd{}{\tilde\phi} : a=1,\dots,n-1\Big\}
\end{equation*}
This is a $2$-step nilpotent Lie algebra of complex dimension $2n$, whose center is 
\begin{equation*}
	\mf{z}(\mf{h}) = \Span_\C\Big\{V_a, \pd{}{\tilde\phi} : a=1,\dots,n-1\Big\}.
\end{equation*}
The only non-trivial brackets are $[Y_a,V_0] = V_a$, for $a=1,\dots,n-1$. Analogous statements hold for $\overline{\mf{h}}$. Note that $\mf{h}\cap\overline{\mf{h}} = \C\pd{}{\tilde\phi}$, and that $\mf{h} + \bar{\mf{h}}$ generates the Lie algebra $\mf{g}^\C$. Note that the elements of $\mf g^\C$ preserve the coordinate function $\rho$ and are therefore tangent to its level sets, which we shall denote by $\bar N_\rho$.

\begin{lem}\label{framehbarh:lem}
	There is a natural isomorphism $T^\C \bar N_\rho \cong \bar N_\rho \times (\mf{h}+\bar{\mf{h}})$. 
\end{lem}
\begin{myproof}
	We show that the $4n-1$ vector fields $Y_a,V_a,V_0,\bar Y_a,\bar V_a,\bar V_0, \pd{}{\tilde\phi}$ give a global frame for $T^\C \bar N_\rho$. Comparing dimensions, it is enough to check that these complex vector fields are pointwise linearly independent. We have the decomposition
	\begin{equation*} 
		T^\C \bar N_\rho \cong (\mathrm{pr}_X^*TB^{n-1})^\C \oplus (\mathrm{pr}_w^*T\C^n)^\C 
		\oplus \C\pd{}{\tilde\phi},
	\end{equation*}
	where $B^{n-1}\subset \C^{n-1}$ denotes the open unit ball, and we used the natural projections $\mathrm{pr}_X: \bar N_\rho = B^{n-1}\times\C^n\times\R\to B^{n-1}$, $\mathrm{pr}_w: \bar N_\rho = B^{n-1}\times\C^n\times\R\to \C^{n}$.
	Now consider a linear combination 
	\begin{equation*}
		0 = \sum_{a=1}^{n-1}\lambda_aY_a + \mu_a\bar Y_a + \sum_{k=0}^{n-1}\delta_kV_k + \epsilon_k\bar V_k + \gamma \pd{}{\tilde\phi},
	\end{equation*}
	with coefficients $\lambda_a,\mu_a,\delta_k,\epsilon_k,\gamma\in\C$. Projecting onto the $(TB^{n-1})^\C = T^{1,0}B^{n-1} \oplus T^{0,1}B^{n-1}$ component and considering $(1,0)$ and $(0,1)$ parts separately yields the relations 
	\begin{equation*}
		\mu_a = X^a\sum_{j=1}^{n-1}X^j\lambda_j,\qquad 
		\lambda_a = \bar X^a\sum_{j=1}^{n-1}\bar X^j\mu_j.
	\end{equation*}
	Using vector notation $\lambda = (\lambda_1,\dots,\lambda_{n-1})\in\C^{n-1}$ and similarly for $\mu$ and $X$, we can use the standard Hermitian inner product on $\C^{n-1}$ to rewrite these relations as $\mu = \langle \lambda,\bar X\rangle X$ and $\lambda = \langle \mu,X\rangle \bar X$. The Cauchy--Schwarz inequality now yields
	\begin{equation*}
		\norm{\lambda}\leq \norm{\mu}\norm{X}^2, \qquad \norm{\mu}\leq\norm{\lambda}\norm{X}^2,
	\end{equation*}
	which implies $(1-\norm{X}^4)\norm{\lambda} =0$ since $\norm{X}<1$ (cf.~the discussion preceding \eqref{eq:fullmetric}). Since $\norm{X}<1$, this means that $\lambda=0$, hence also $\mu=0$. Projecting onto the $T\C^{n-1}$-factor and considering $(1,0)$ and $(0,1)$ parts, we see immediately from \eqref{eq:DefVa} that $\delta_k = 0= \epsilon_k$ for all $k=0,\dots, n-1$, which in turn implies $\gamma=0$, completing the proof. 
\end{myproof}

Consider once again the submanifold $H = \{X = 0\}\subset \bar N_{\rho_0}$ and note that, at $p = (0,w,\tilde\phi,\rho_0)\in H$, we have
\begin{equation*}
	\qkvf_a(p)=\pd{}{\bar X^a} - w^0\pd{}{w^a}-\bar w^a\pd{}{\bar w^0},
\end{equation*}
and 
\begin{equation*} 
	\qkheis_a(p)= \pd{}{w^a} - i\bar w^a\pd{}{\tilde\phi}, \qquad 
	\qkheis_0(p)= \pd{}{w^0} +i\bar w^0 \pd{}{\tilde\phi}. 
\end{equation*}
The metric (given in \eqref{eq:fullmetric}) also simplifies on $T_p\bar N_{\rho_0}$: 
\begin{equation}\label{eq:gcH}
	\begin{aligned} 
		g^c &= \frac{\rho_0 + c}{\rho_0}\sum_{a=1}^{n-1}\big|\d X^a\big|^2 + \frac{1}{2\rho_0}\sum_{a=1}^{n-1}\big|\d w^a\big|^2 + \bigg(\frac{\rho_0+c}{\rho_0^2}-\frac{1}{2\rho_0}\bigg)\big|\d w^0\big|^2\\
		& \quad+ \frac{1}{4\rho_0^2}\frac{\rho_0+c}{\rho_0+2c}\bigg(\d\tilde\phi -2i\bigg(\bar w^0\d w^0 - w^0\d\bar w^0 -\sum_{a=1}^{n-1}(\bar w^a\d w^a-w^a\d\bar w^a)\bigg)\bigg)^2.
	\end{aligned}
\end{equation}

\begin{lem}\label{lem:parallelnormalTH}
	The vector fields $\pd{}{X^a}, \pd{}{\bar X^a}$, $a=1,\dots,n-1$ form a parallel frame for the complexified normal bundle of the submanifold $H\subset \bar N_{\rho_0}$ with respect to the connection induced by the Levi-Civita connection of $(\bar N_{\rho_0}, g^c)$. 
\end{lem}
\begin{myproof}
	By looking at the form of the metric $g^c$ at $p\in H$ given in \eqref{eq:gcH} we see that the vector fields $\pd{}{X^a}, \pd{}{\bar X^a}$, $a=1,\dots,n-1$ are orthogonal to $T^\C_pH$. Moreover, the complex bilinear extension of the metric satisfies
	\begin{align*}
		g^c\bigg(\pd{}{X^a},\pd{}{X^b}\bigg) &= 0 \\
		g^c\bigg(\pd{}{\bar X^a},\pd{}{\bar X^b}\bigg) &= 0 \\
		g^c\bigg(\pd{}{X^a},\pd{}{\bar X^b}\bigg) &= \frac{\rho_0 +c}{\rho_0}\delta_{ab}.
	\end{align*}
	Now it follows from the Koszul formula that for any vector field $W\in\mf{X}^\C(H)$ the covariant derivatives $\nabla_{W} \pd{}{X^b}$ and $\nabla_{W} \pd{}{\bar X^b}$ must be tangent to $H$. 
\end{myproof}

Given a vector field $K$ on $\bar N_{\rho_0}$, we shall denote the normal component of the restriction of $K$ to the submanifold $H$ by $K_H^\perp$. 
	
\begin{lem}
	Denote by $\mf{k} = \{K \in \mf X^\C(\bar N_{\rho_0}) : \nabla^\perp(K_H^\perp) = 0 \}$, the space of complex vector fields on $\bar N_{\rho_0}$ whose normal component along the submanifold $H$ is a parallel section of the normal bundle of $H$. Then $\mf{g}\subset \mf{k}$.
\end{lem}
\begin{myproof}
	We can write any complex vector field in the form
	\begin{equation*}
		K = \sum_a f_a\pd{}{X^a} + g_a\pd{}{\bar X^a} + \sum_k\phi_k\pd{}{w^k} + 
		\psi_k\pd{}{\bar w^k} + h\pd{}{\tilde\phi}
	\end{equation*}
	with $f_a,g_a,\phi_a,\psi_k,h\in C^\infty(\bar N_{\rho_0})$. By \Cref{lem:parallelnormalTH}, we see that $K\in\mf{k}$ if and only if the functions $f_a,g_a$ are constant when restricted to the submanifold $H$. Clearly, this is satisfied by the vector fields $\qkvf_C, \qkvf_a,\bar \qkvf_a, \qkheis_k,\bar \qkheis_k,\pd{}{\tilde\phi}$. Moreover, at $p = (0,w,\tilde\phi,{\rho_0})$ we have
	\begin{equation*}
		[\qkvf_a,\bar \qkvf_b](p) = \delta_{ab}\Bigg(  w^0\pd{}{w^0}- \bar w^0\pd{}{\bar w^0}
		-2ic\pd{}{\tilde\phi}\Bigg)  + \bar w^a \pd{}{\bar w^b}  - w^b\pd{}{w^a},
	\end{equation*}
	which is also tangent to $H$ and thus is an element of $\mf{k}$. This means that we have verified our claim for the basis $\{\qkvf_C, \qkvf_a,\bar{\qkvf_a},[\qkvf_a,\bar \qkvf_b],\qkheis_k,\bar{\qkheis_k},\pd{}{\tilde\phi}\}$ of $\mf{g}^\C$.
\end{myproof}

\begin{prop}\label{prop:char}
	For every $c>0$, the elements of the algebra $\mf{g}$ can be characterized, modulo Killing fields that are zero at each point of $H$, as those Killing fields on $(\bar N_{\rho_0},g^c)$ whose normal component along $H$ is parallel. 
\end{prop}
\begin{myproof}
	Let $K$ be such a Killing field. Expand it in the frame $T^\C \bar N_{\rho_0} \cong \bar N_{\rho_0} \times (\mf{h}+\bar{\mf{h}})$ constructed in Lemma~\ref{framehbarh:lem}. Then the coefficient functions in front of $Y_a,\bar Y_a$ are constant along $H$. By subtracting a linear combination of $Y_a,\bar Y_a$ with constant coefficients, we can thus obtain a Killing field that is tangent to $H$. By Lemma~\ref{tgH:lem}, this must then coincide with an element of $\alpha((\mf{u}(1)\times\mf{u}(n-1))\ltimes\heis_{2n+1})$ along $H$. This shows that, up to  addition of a Killing field vanishing along $H$, the vector field $K$ lies in $\mf{g}$.
\end{myproof}

\subsection{The action of \texorpdfstring{$\widetilde{\U}(1,n-1)\ltimes \Heis_{2n+1}$}{the group}}
\label{sec:isom}

In the previous section we have constructed an injective ($c$-dependent) infinitesimal action $\alpha: \mf{u}(1,n-1)\ltimes\mf{heis}_{2n+1}\to\mf{g}\subset \mf{isom}(\bar N_n,g^c)$. Since $(\bar N, g^c)$ is a complete Riemannian manifold, any Killing field on it is necessarily complete. We can therefore integrate the infinitesimal action $\alpha: \mf{u}(1,n-1)\ltimes\heis_{2n+1}\to \mf{isom}(\bar N_n,g^c)$ to an isometric action of the corresponding simply connected group $\widetilde{\U}(1,n-1)\ltimes \Heis_{2n+1}$. In this way we obtain a group homomorphism $\beta:\widetilde{\U}(1,n-1)\ltimes \Heis_{2n+1}\to\Isom(\bar N_n,g^c)$. Let us collect some useful pieces of notation, which we will repeatedly make use of in the following.

\begin{not*}
	We will abbreviate $\tilde{\mc G}\coloneqq\widetilde{\U}(1,n-1)\ltimes \Heis_{2n+1}$ and denote its image $\beta(\tilde{\mc G})$ by $\mc G$. The subgroups $\beta(\widetilde{\U}(1,n-1))$ and $\beta(\Heis_{2n+1})$ will be called $\mc U$ and $\mc H$, respectively. We will write $\mc F$ for the intersection $\mc U\cap \mc H$. We will also consider $\tilde{\mc G}'\coloneqq \widetilde{\SU}(1,n-1)\ltimes\Heis_{2n+1}$; the restriction of $\beta$ to $\tilde{\mc G'}$ will be denoted by $\beta'$. Its image is $\mc G'\coloneqq\beta'(\tilde{\mc G}')$, and we abbreviate $\beta'(\widetilde{\SU}(1,n-1))$ to $\mc U'$ (note that $\beta'(\Heis_{2n+1})=\mc H$). Finally, we set $\mc F'\coloneqq \mc U'\cap \mc H$.
\end{not*}

We are interested in $\mc G$, which acts effectively on $\Isom(\bar N_n,g^c)$. In order to understand $\mc G$, we have to determine the kernel of $\beta$.
 
Since the infinitesimal action $\alpha:\mf u(1,n-1)\ltimes \heis_{2n+1}\to \isom(\bar N_n,g^c)$ is injective, $\ker\beta$ is discrete and normal, hence central. The elements in $\ker\beta$ stabilize every point of $\bar N_n$, and in particular the point $n_0=(0,0,0,\rho_0)$, for any fixed $\rho_0>0$. Thus, $\ker\beta$ must be a discrete central subgroup of the stabilizer $\tilde{\mc{G}}_{n_0}$. We computed its Lie algebra $\mf{g}_{n_0}$ in \Cref{lem:stabilizer}, and it turns out that we can integrate the vector fields in the center of $\mf{g}_{n_0}$ explicitly. This allows us to determine $\ker\beta$, as we will now show.

\begin{lem}\label{lem:kerbeta}
	For $n\geq 2$, the kernel of $\beta$ is the subgroup of $Z(\widetilde{\U}(1,n-1)\ltimes\Heis_{2n+1}) \cong \R\times2\pi\Z\times \R$ generated by $(2\pi,0,4\pi c)$ and $\big(-\frac{2\pi}{n},-2\pi,4\pi\frac{n-2}{n}c\big)$. In particular, $\ker\beta\cong \Z^2$.
\end{lem}
\begin{myproof}
	As shown in \Cref{lem:stabilizer}, the Lie algebra $\mf g_{n_0}$ of the stabilizer of $n_0$ is isomorphic to $\mf u(1)\times \mf u(n-1)$. Its two-dimensional center is spanned by the vector fields 
	\begin{equation}\label{eq:centralvectors}
		\begin{aligned}
			\mc C_1&\coloneqq \qkvf_C+2c\pd{}{\tilde\phi},\\
			\mc C_2&\coloneqq \sum_{a=1}^{n-1} \Im([Y_a,\bar Y_a])+2(n-1)c\pd{}{\tilde\phi}-\mc C_1,
		\end{aligned}
	\end{equation}
	as can be easily checked using \eqref{eq:LieUaUbsigma}. Let us reproduce the explicit coordinate expressions, obtained from \eqref{eq:alphagln} and \eqref{eq:commexpression}, for convenience:
	\begin{align*}
		\mc C_1&=-i\sum_{k=0}^{n-1}\bigg(w^k\pd{}{w^k}-\bar w^k\pd{}{\bar w^k}\bigg),\\
		\mc C_2&=-in\sum_{a=1}^{n-1}\bigg(X^a\pd{}{X^a}-\bar X^a\pd{}{\bar X^a}+w^0\pd{}{w^0}-\bar w^0\pd{}{\bar w^0}\bigg).
	\end{align*}
	Integrating these vector fields, we see that the corresponding one-parameter groups of diffeomorphisms act as follows:
	\begin{align*}
		\Phi_t^{\mc C_1}(X^a,w^0,w^a,\tilde\phi,\rho)&=(X^a,e^{-it}w^0,e^{-it}w^a,\tilde\phi,\rho)\\
		\Phi_t^{\mc C_2} (X^a,w^0,w^a,\tilde\phi,\rho)&=(e^{-int}X^a,e^{-int}w^0,w^a,\tilde\phi,\rho)
	\end{align*}
	It is clear that the flows are periodic with periods $2\pi$ and $\frac{2\pi}{n}$, respectively. Thus, the $\Z^2$-subgroup generated by $\exp(2\pi \mc C_1)$ and $\exp(\frac{2\pi}{n}\mc C_2)$ acts trivially on all of $\bar N_n$ and is therefore contained in $\ker \beta$. In fact, it is equal to $\ker\beta$, since we know that $\ker \beta$ is certainly contained in the preimage of the stabilizer.
	
	The center of $\widetilde{\U}(1,n-1)\ltimes \Heis_{2n+1}$ is
	\begin{equation*}
		Z(\widetilde{\U}(1,n-1)\ltimes\Heis_{2n+1})
		=\exp(\R\cdot C)\times\exp(2\pi \Z\cdot C')\times \exp(\R\cdot T)
		\cong \R\times2\pi \Z\times\R,
	\end{equation*}
	where $C'\coloneqq \frac{i}{n} \operatorname{diag}(1-n,1,\dots,1)\in \mf {su}(1,n-1)$. 
	
	Using \eqref{eq:centralvectors} and combining the fact that $\sum_a \Im_\sigma [U_a,U_a^\sigma]=n C'$ and that $\alpha$ is an anti-homomorphism, we can directly read off that the kernel of $\beta$ is the subgroup generated by $(2\pi,0,4\pi c)$ and $\big(-\frac{2\pi}{n},-2\pi,4\pi\frac{n-2}{n}c\big)$.
\end{myproof}

\begin{rem}
	We observe that, for all $c>0$, the subgroup $\ker\beta\subset Z(\tilde{\mc G})$ has a non-trivial projection onto the subgroup $Z(\Heis_{2n+1})$.
\end{rem}

We will now consider the subgroups $\mc U $ and $\mc H$ of $\mc G$. 

\begin{lem}\label{lem:betaUcapH}
	For every $n\geq 2$, the following hold:
	\begin{numberedlist}
		\item The restricted map $\beta: \Heis_{2n+1}\to \mc H$ is an isomorphism. 
		 \item The kernel $\ker\beta$ is contained in the subgroup $\widetilde{\U}(1,n-1)$ if $c=0$. If $c>0$, then the kernel of the restricted map $\beta: \widetilde{\U}(1,n-1)\to \mc U$ is equal to the infinite cyclic group $2\pi\Z\cdot (n-1,n) \subset \R\times \Z = Z(\widetilde{\U}(1,n-1))$ if $n$ is odd and $\pi\Z\cdot (n-1,n)$ if $n$ is even.
		\item The intersection $\mc F=\mc U\cap \mc H$ is given by $\langle \beta(0,0,\frac{8\pi c}{n})\rangle$ if $n$ is even and $\langle \beta(0,0,\frac{4\pi c}{n})\rangle$ if $n$ is odd. In particular, it is trivial if $c=0$ and infinite cyclic if $c>0$.   
	\end{numberedlist}
\end{lem}
\begin{myproof}\leavevmode
	\begin{numberedlist}
		\item By \Cref{lem:kerbeta}, the kernel of $\beta$ is the subgroup of $Z(\tilde{\mc G}) \cong \R\times2\pi\Z\times \R$ generated by $(2\pi,0,4\pi c)$ and $\big(-\frac{2\pi}{n},-2\pi, 4\pi\frac{n-2}{n}c\big)$. But this group has trivial intersection with $\Heis_{2n+1}$ for any value of $c$, so the restriction of $\beta$ is injective, hence an isomorphism.
	    \item We compute $\ker\beta\cap Z(\widetilde{\U}(1,n-1))$, using the generators provided by \Cref{lem:kerbeta}. For $c=0$, our claim is obvious. Now assume that $c>0$ and let $x,y\in \Z$. The third component of $x(2\pi, 0, 4\pi c) + y\big(-\frac{2\pi}{n}, -2\pi, 4\pi\frac{n-2}{n}c\big)$ is $4\pi c(x+\frac{n-2}{n}y)$, which vanishes if and only if $nx + (n-2)y = 0$, which in turn holds if and only if $(x,y)$ is an integral multiple of $(2-n,n)$ if $n$ is odd and $\frac{1}{2}(2-n,n)$ if $n$ is even. Thus, if $n$ is odd, then  $\ker\beta$ is generated by the element  $(n-2)\cdot(2\pi,0,4\pi c) - n\cdot\big(-\frac{2\pi}{n},-2\pi,4\pi\frac{n-2}{n}c\big) = 2\pi\cdot(n-1,n,0)$. If $n$ is even we obtain $\pi\cdot(n-1,n,0)$ as a generator instead. 
		\item Since $\dim\mc G = \dim \mc U +\dim \mc H$, $\mc F$ is discrete and normal in $\mc U$. Thus 
		$\mc U \subset Z(\mc U)$. Moreover, since left-translation by $\mc U$ fixes the origin in $\R^{2n}\times\{0\}\subset \mc H$, we must have $\mc F\subset Z(\mc H)$. It follows that $\mc F = Z(\mc U)\cap Z(\mc H)$. The preimage $\beta^{-1}(\mc F) = \beta^{-1}(Z(\mc U))\cap\beta^{-1}(Z(\mc H))$ is therefore a discrete subgroup of the center $\R\times 2\pi\Z \times\R$ of $\tilde{\mc{G}}$.
				
		By part (i), $\beta$ restricts to an isomorphism $\Heis_{2n+1}\to \mc H$. Consequently, we have $\beta^{-1}(Z(\mc H))=(\{(0,0)\}\times \R)+\ker \beta$. On the other hand $\beta^{-1}(Z(\mc U)) = (\R\times 2\pi\Z\times\{0\}) + \ker\beta$. Now, a point $2\pi\cdot (0,0,\lambda)$ with $\lambda\in\R$ lies in $\beta^{-1}(\mc F)$ if and only if we can find some $\kappa\in\ker\beta$, $a\in \R$ and $b\in \Z$ such that $2\pi\cdot (0,0,\lambda) + \kappa= 2\pi\cdot (a,b,0)$. 
		For any $x,y\in \Z$ we have
		\begin{align*}
			(0,0,2\pi\lambda)\! + x(2\pi,0,4\pi c)\! - y(\tfrac{2\pi}{n},2\pi,-4\pi\tfrac{n-2}{n}c)
			= 2\pi \big(x\!-\!\tfrac{y}{n},-y,\lambda\! +\! 2c(x+\tfrac{n-2}{n}y)\big).
		\end{align*}
		To make the third entry zero, we must have  $\lambda = -2c\big(x+\frac{n-2}{n}y\big)\in 2c\cdot\Span_{\Z}\{1, \frac{n-2}{n}\}$, which equals $\frac{4c}{n}\Z$ if $n$ is even and $\frac{2c}{n}\Z$ if $n$ is odd. Thus, $\beta^{-1}(\mc F) = \Z\cdot(0,0,\frac{8\pi c}{n}) + \ker\beta$ if $n$ is even and $\beta^{-1}(\mc F)=\Z\cdot(0,0,\frac{4\pi c}{n})+\ker\beta$ if $n$ is odd, which means that $\mc F=\langle\beta(0,0,\frac{8\pi c}{n})\rangle$ in the former case and $\mc F=\langle\beta(0,0,\frac{4\pi c}{n})\rangle$ in the latter, just as claimed.\qedhere
	\end{numberedlist}
\end{myproof}

Next, we study $\mc G'$, which is obtained by restricting $\beta$ to $\tilde{\mc G}'=\widetilde{\SU}(1,n-1)\ltimes\Heis_{2n+1}$.

\begin{lem}\label{lem:betaU'capH}
	For every $n\geq 2$, the following hold:
	\begin{numberedlist}
		\item If $c=0$, then the kernel of $\beta':\widetilde{\SU}(1,n-1)\to \mc U'$ is equal to $2\pi n \Z\subset 2\pi \Z = Z(\widetilde{\SU}(1,n-1))$. Thus, $\mc U'\cong \SU(1,n-1)$.
		\item If $c>0$ then $\beta': \widetilde{\SU}(1,n-1)\to \mc U'$ is an isomorphism. 
		\item The subgroup $\mc F'=\mc U'\cap \mc H$ is given by $\langle \beta'(0,0,4\pi c(n-1))\rangle$, hence trivial if $c=0$ and infinite cyclic if $c>0$.
	\end{numberedlist}
\end{lem}
\begin{myproof}\leavevmode
	\begin{numberedlist}
		\item Since $c=0$, \Cref{lem:kerbeta} implies that the  kernel of $\beta$ intersects $Z(\widetilde{\SU}(1,n-1)) \cong 2\pi\Z$ in the cyclic subgroup generated by $2\pi n$. The assertion follows.
		\item We have $\widetilde{\U}(1,n-1) \cong \R\times\widetilde{\SU}(1,n-1)$ and $\mc U' = \beta(\{0\}\times\widetilde{\SU}(1,n-1))$. By part (ii) of \Cref{lem:betaUcapH}, the intersection of $\ker\beta$ and $\{0\}\times \widetilde{\SU}(1,n-1)$ is trivial.
		\item By the same arguments as in the proof of part (iii) of \Cref{lem:betaUcapH}, $\mc F' = Z(\mc U')\cap Z(\mc H)$. The center of $\tilde{\mc G}' = (\{0\}\times\widetilde{\SU}(1,n-1))\ltimes \Heis_{2n+1}$ is $\{0\}\times 2\pi\Z\times\R$, which intersects the kernel of $\beta$ in the cyclic subgroup $\Z\cdot(0,-2\pi n,4\pi c(n-1))$. This observation together with a calculation analogous to the proof of part (iii) of \Cref{lem:betaUcapH} shows that $\mc F' = \langle \beta'(0,0,4\pi c(n-1))\rangle$.\qedhere
\end{numberedlist}
\end{myproof}

The following Proposition is now an immediate application of the results proven in \Cref{lem:kerbeta,lem:betaUcapH,lem:betaU'capH}.

\begin{prop}\label{prop:G'structure}
	Let $n\geq 2$ and consider $\mc G'\subset \Isom(\bar N_n,g^c)$.
	\begin{numberedlist}
		\item If $c=0$, then $\mc G' \cong \SU(1,n-1)\ltimes\Heis_{2n+1}$. 
		\item If $c>0$, then $\mc G' \cong (\widetilde{\SU}(1,n-1)\ltimes \Heis_{2n+1})/\mc F'$,
		where $\mc F'$ is the infinite cyclic subgroup $\langle \beta'(0,0,4\pi c(n-1))\rangle$.
		\proofclear
	\end{numberedlist}
\end{prop}

For $c>0$ and $n\geq 2$ we now consider the cyclic quotient $\hat N_n^c\coloneqq \bar N_n/\mc F'$, endowed with the induced quaternionic K\"ahler metric which we continue to denote by $g^c$. This amounts to considering the coordinate $\tilde\phi$ periodic with period $4\pi c(n-1)$. We extend this family to $c=0$ by setting $\hat N_n^0=\bar N_n$. 

Since $\mc F'\subset \tilde{\mc G}$ is central, the action of $\tilde{\mc G}$ on $\bar N_n$ descends to $\hat N_n^c$. We extend our notation to the resulting action in the obvious way:

\begin{not*}
	We denote by $\hat\beta:\tilde{\mc G}\to \Isom(\hat N_n^c,g^c)$ and $\hat\beta': \tilde{\mc G}'\to \Isom(\hat N_n^c, g^c)$ the corresponding group homomorphisms with images $\hat{\mc G}$ and $\hat{\mc G}'$, respectively. We also define $\hat{\mc H}\coloneqq \hat\beta(\Heis_{2n+1})$ and $\hat{\mc U}=\hat\beta(\widetilde{\U}(1,n-1))$, as well as $\hat{\mc U}'=\hat\beta'(\widetilde{\SU}(1,n-1))$. 
\end{not*}

We can now state the main result of this section.

\begin{thm}\label{thm:SUltimesHaction}
	For $n\geq 2$ and any value of $c\geq 0$, $\hat{\mc G}'\subset \Isom(\hat N_n^c,g^c)$ is isomorphic to $\SU(1,n-1)\ltimes(\Heis_{2n+1}/\mc F')$. In particular, the group $\SU(1,n-1) \ltimes (\Heis_{2n+1}/\mc F')$ acts effectively and isometrically on $(\hat N_n^c,g^c)$.
\end{thm}
\begin{myproof}
	By construction, and the proof of part (iii) of \Cref{lem:betaU'capH}, the kernel of $\hat \beta'$ is generated by $(0,2\pi n,0)$ and $(0,0,4\pi c(n-1))$. This means that $\im\hat\beta'\cong (\widetilde{\SU}(1,n-1)\ltimes\Heis_{2n+1})/\ker\hat\beta'\cong \SU(1,n-1)\ltimes(\Heis_{2n+1}/\mc F')$, as claimed.
\end{myproof}

With the notation introduced above, we may rephrase this result as follows. Firstly, $\hat{\mc U}'\cong \SU(1,n-1)$ and $\hat{\mc H}\cong \Heis_{2n+1}/\mc F'$. Secondly $\hat{\mc U}'\cap \hat{\mc H}=\{\id\}$ so that $\hat{\mc G}'\cong \hat{\mc U}'\ltimes \hat{\mc H}$.

The above results have been formulated for the case $n\geq 2$ only. The arguments in the case $n=1$ are very similar and the differences are mostly notational. We state them separately in the next Proposition. 

\begin{prop}\label{prop:betaUcapHn=1}
	Let $n=1$ and $c\geq 0$ arbitrary. Then
	\begin{numberedlist}
		\item $\ker\beta = \Z\cdot  (2\pi, 4\pi c)\subset \R\times\R \cong Z(\widetilde{\U(1)}\ltimes \Heis_3)$.
		\item The restricted map $\beta: \Heis_3\to \mc H$ is an isomorphism.
		\item The kernel of the restricted map $\beta: \widetilde{\U(1)}\to \mc U$ is $2\pi\Z\subset\R = \widetilde{\U(1)}$ if $c =0$ and trivial if $c>0$.
		\item The intersection $\mc F = \mc U\cap \mc H$ is given by $\langle\beta(0, 4\pi c)\rangle$. In particular, it is trivial if $c=0$ and infinite cyclic if $c>0$.
		\item $\mc G = (\widetilde{\U(1)}\ltimes\Heis_{3})/\mathcal F$.
		\item $\mc G' = \Heis_3$.
		\proofclear
	\end{numberedlist}
\end{prop}

In the following, it will be convenient to set $\hat N_1^c\coloneqq \bar N_1/\mc F$ for $c>0$ and $\hat N_1^0\coloneqq\bar N_1$ so that we can speak of $\hat N_n^c$ for all $c\geq 0$ and $n\in\N$.

Our next aim is to prove that $\mc G$ and $\hat{\mc G}$ are closed subgroups of $\Isom(\bar N_n,g^c)$ and $\Isom(\hat N_n^c,\hat g^c)$, respectively. We first consider the subgroup $\hat{\mc U}=\hat\beta(\widetilde{\U}(1,n-1))\subset \hat{\mc G}$.

\begin{lem}
	For all $n\in \N$ and any $c\geq 0$, $\hat{\mc U}\subset \Isom(\hat N_n^c,g^c)$ has compact center. 
\end{lem}
\begin{myproof}
	The $\mf u(1)$-factor of $\Lie(\hat{\mc U})\cong \mf u(1)\oplus\mf{su}(1,n-1)$ is generated by the vector field $Y_C$. If $c=0$, it is clear that $Y_C$ generates a $\U(1)$-subgroup. In the case $c>0$, the periodicity of the coordinate $\tilde\phi$ ensures that this assertion remains valid, so $\hat{\mc U}=\U(1)\cdot\hat{\mc U}'$, where $\hat{\mc U}'$ is the group of isometries generated by the subalgebra $\mf{su}(1,n-1)$. Since the $\U(1)$-factor is central, we have $Z(\hat{\mc U})=\U(1)\cdot Z(\hat{\mc U}')$, and it therefore suffices to prove that $Z(\hat{\mc U}')$ is compact. In the case $n=1$, $\hat{\mc U}'$ is trivial so there is nothing to prove. For $n\geq 2$, \Cref{thm:SUltimesHaction} shows that $\hat{\mc U}'\cong \SU(1,n-1)$, so its center is cyclic of order $n$ and in particular finite.  
\end{myproof}

For all $c>0$, the periodicity of the coordinate $\tilde\phi$ implies that the center of $\hat{\mc H}=\hat\beta(\Heis_{2n+1})$ is compact as well. By~\cite[Prop.~4.2]{Mac1993} the compactness of the centers of $\hat{\mc U}$ and $\hat{\mc H}$ implies that they are closed subgroups of $\Isom(\hat N_n^c,g^c)$. 

\begin{prop}
	 For all $n\in \N$ and $c\geq 0$, $\hat{\mc G}$ is a closed subgroup of $\Isom(\hat N_n^c,g^c)$ and $\mc G$ is a closed subgroup of $\Isom(\bar N_n,g^c)$.
\end{prop}
\begin{myproof}
	We start by proving the result in the case $c>0$. The idea is to exploit the closedness of the subgroups $\hat{\mc U}$ and $\hat{\mc H}$ of $\hat{\mc G}$ to prove the first claim, and then apply a lifting argument to establish the second claim.
	
	Note that $\hat{\mc H}\subset \hat{\mc G}$ is normal, i.e.~$\hat{\mc G}$ is contained in the normalizer $N(\hat{\mc H})$,  which is a closed subgroup, since it coincides with the normalizer of the Lie algebra of $\hat{\mc H}$. Thus, it suffices to prove that $\hat{\mc G}$ is closed as a subgroup of $N(\hat{\mc H})$. Let us therefore consider a sequence $\{x_j\}$ in $\hat{\mc G}$ which converges in $N(\hat{\mc H})$. Under the projection $N(\hat{\mc H})\to N(\hat{\mc H})/\hat{\mc H}$, which is a submersion because $\hat{\mc H}$ is closed, this projects to a sequence in $\hat{\mc U}/(\hat{\mc U}\cap \hat{\mc H})\subset N(\hat{\mc H})/\hat{\mc H}$. 
	
	It follows from the description of $\mc U\cap \mc H$ given in \Cref{lem:betaUcapH} and \Cref{prop:betaUcapHn=1}, and the definition of $\hat N_n^c$, that $\hat{\mc U}\cap \hat{\mc H}$ is a finite cyclic subgroup, so the map $\hat{\mc U}\to \hat{\mc U}/(\hat{\mc U}\cap \hat{\mc H})$ is a finite cyclic covering. The center of $\hat{\mc U}/(\hat{\mc U}\cap \hat{\mc H})$ is compact (because the center of $\hat{\mc U}$ is) and we may appeal once again to \cite[Prop.~4.2]{Mac1993} to deduce that it is a closed subgroup of $N(\hat{\mc H})/\hat{\mc H}$. The projected sequence therefore converges in $\hat{\mc U}/(\hat{\mc U}\cap \hat{\mc H})$, and lifts to the covering space, producing a convergent sequence $\{g_j\}$ in $\hat{\mc U}$. We may now write $x_j=g_jh_j$, and it follows from convergence of $\{x_j\}$ and $\{g_j\}$ that $\{h_j\}$ converges in $N(\hat{\mc H})$. In fact, since $\hat{\mc H}$ is closed in $N(\hat{\mc H})$ it converges in $\hat{\mc H}$, and we conclude that $\{x_j\}$ converges in $\hat{\mc G}=\hat{\mc U}\cdot\hat{\mc H}$. This proves that $\hat{\mc G}\subset \Isom(\hat N_n^c,\hat g^c)$ is closed.
	
	By definition, $\hat{\mc G}=\mc G/\mc F'$ where $\mc F'\cong\Z$ is a discrete subgroup of the center $\mc Z$ of $\mc H=\beta(\Heis_{2n+1})$. Let $C(\mc Z)\subset \Isom(\bar N_n,g^c)$ denote its centralizer. Then the covering map $\bar N_n\to \hat N_n^c$ induces a covering $C(\mc Z)\to C(\mc Z/\mc F')$. Under this covering map, $\mc G$ is the preimage of the closed subgroup $\hat{\mc G}\subset C(\mc Z/\mc F')$, and therefore $\mc G$ is itself closed in $C(\mc Z)$, hence in $\Isom(\bar N_n,g^c)$ as well.
	
	In the case $c=0$ we proceed analogously, by passing to the quotient $(\bar N_n/\Z,g^0)$ obtained by making $\tilde\phi$ periodic with period $2\pi$, showing that $\mathcal G/\Z$ is a closed subgroup of $\Isom(\bar N_n/\Z, g^0)$ and then performing the same lifting argument.
\end{myproof}

We note the following consequence, which we will use later.

\begin{cor}\label{cor:compactstab}
	Let $n\in\N$ and $c\geq 0$ and consider arbitrary points $p\in \bar N_n$ and $\hat p\in \hat N_n^c$. Then $\mc G_0\coloneqq \{g\in \mc G\mid g\cdot p=p\}$ and $\hat{\mc G}_0\coloneqq \{\hat g\in \hat{\mc G}\mid \hat g\cdot \hat p =\hat p\}$ are compact.\proofclear
\end{cor}

\section{Quotients and ends of finite volume}
\label{sec:quotients}

\subsection{Quotients from arithmetic lattices}

We will now use the group actions discussed in the previous section to construct interesting quotients of $(\bar N_n,g^c)$ by dividing out appropriate discrete subgroups of isometries. The resulting manifolds will be complete and quaternionic K\"ahler (since the isometries automatically preserve the quaternionic structure), with non-trivial fundamental group.

For $c=0$, the metric is symmetric and much is known about its quotients by discrete subgroups. In particular, it is known that $(\bar N_n,g^0) =\frac{\SU(n,2)}{\mathrm{S}(\U(n)\times \U(2))}$ admits a compact quotient \cite{Bor1963}. Since any isometry of $(\bar N_n,g^c)$, $c>0$, preserves% 
\footnote{This follows from the fact that $\rho$ is a curvature invariant for $c>0$ \cite{CST2022}. For $c=0$, this is of course not true, as witnessed by the fact that $(\bar N_n,g^c)$ is the symmetric space $\frac{\SU(n,2)}{\mathrm{S}(\U(n)\times\U(2))}$.} %
the level sets of the global coordinate function $\rho$, we cannot expect to obtain compact quotients in this case. The best one can hope for is to obtain manifolds that have the structure of a fiber bundle over $\R_{>0}$, which is parametrized by $\rho$, with locally homogeneous and compact fibers. 

If $n\leq 2$ we shall explicitly construct (infinitely many) discrete subgroups of $\Isom(\bar N_n,g^c)$, $c>0$, which yield complete quaternionic K\"ahler manifolds of this type. We shall see that if one does not insist on compactness of the fibers and only requires them to be of finite volume, suitable discrete subgroups can be found for all values of $n\in \N$.

For the purpose of studying quotients of $(\bar N_n,g^c)$, we may (recalling the notation of \Cref{sec:isom}) restrict our attention to the subgroup $\mc G'\subset \mc G$, whose orbits are the same as those of the full group. This is most easily seen by noting that $\mf g=\Lie(\mc G)=\mf u(1)\oplus \Lie(\mc G')$, where the first factor is generated by the vector field $Y_C$, which is tangent to the fibers of the projection $(\bar N_n,g^c)\to \CH^n$ and preserves the coordinate function $\rho$. But the level sets of $\rho$ in these fibers are already acted upon transitively by $\mc H\cong \Heis_{2n+1}$. As we will see shortly, restricting to $\mc G'$ affords us certain technical advantages.

Let us now turn to the construction of discrete subgroups. If $c=0$, then $\mc G'$ has the structure of a semi-direct product (cf.~\Cref{prop:G'structure}). This no longer holds true for $c>0$, but can be remedied by passing to the cyclic quotient $\hat N_n^c$: by \Cref{thm:SUltimesHaction}, the induced group of isometries $\hat{\mc G}'$ is always a semi-direct product of the form $\SU(1,n-1)\ltimes \hat{\mc H}$, where $\hat{\mc H} =  \Heis_{2n+1}/\mc F'$, with $\mc F'\cong \Z$ a central, cyclic subgroup of $\Heis_{2n+1}$. 

Thinking of $\hat{\mc G}'$ as determined by an action of $\SU(1,n-1)$ on $\hat{\mc H}$ by automorphisms suggests the following strategy. First, we construct a discrete subgroup $\bar\Gamma_1$ of $\SU(1,n-1)$. Next, we construct a discrete subgroup $\hat \Gamma_2\subset\hat{\mc H}$ which is invariant under the action of $\bar\Gamma_1$. Then $\hat\Gamma\coloneqq \bar\Gamma_1\ltimes \hat \Gamma_2$ will be a discrete subgroup of the semi-direct product group.

The first step, thus, is to construct a suitable discrete subgroup of $\SU(1,n-1)$. Recall that a lattice in a unimodular Lie group $G$ is a discrete subgroup $\Gamma\subset G$ such that $\vol(G/\Gamma)$ is finite. Here $\vol(G/\Gamma) = \vol(F)$, where $F\subset G$ is a fundamental domain and $\vol(F)$ is computed with respect to the Haar measure on $G$. A lattice is called co-compact if $G/\Gamma$ is compact. The study of lattices in semi-simple Lie groups is a well-developed area of research (see e.g.~\cite{Mor2015}), and we will now show how to apply some of the known constructions to the problem at hand.

In the case $n=1$ the group $\SU(1,n-1)$, and thus the first step in our construction, is trivial. We can take $\hat\Gamma = \hat \Gamma_2$ where $\hat \Gamma_2$ is a discrete subgroup in $\hat{\mc H}=\Heis_3/\mc F'$. If $\hat\Gamma$ is a co-compact lattice in $\hat{\mc H}$, the resulting quotient  $\hat{N}_1/\hat \Gamma$ is topologically $\R_{>0}\times \hat{{\mc H}}/\hat \Gamma$ and hence of the desired type.

The next case, $n=2$, is much less simple. The first step is to construct co-compact lattices in $\SU(1,1)$. It is well-known that $\SU(1,1)$ is isomorphic to $\mathrm{SL}(2,\R)$, so discrete subgroups of $\SU(1,1)$ are the same thing as Fuchsian groups. Therefore, we are looking for co-compact Fuchsian groups which preserve a lattice in $\Heis_5$. One way to construct Fuchsian groups is via quaternion algebras, and this is what we recall next. Let $F$ be a field, and let $a,b\in F^* = F\setminus\{0\}$. Then the quaternion algebra over $F$ associated to the pair $(a,b)$ is the unital algebra over $F$ with generators $I,J,K$ and relations 
\begin{equation*}
	I^2 = a, \quad J^2 = b, IJ = K = -JI
\end{equation*} 
Note that this implies $K^2 = -IJJI = -ab$. We will denote this algebra by $\big(\frac{a,b}{F}\big)$. 
The case relevant for  us is $F = \Q$. It is known that if $a,b$ are positive integers such that $b$ is prime and $a$ is a quadratic non-residue mod $b$, then $\big(\frac{a,b}{\Q}\big)$ is a division algebra (see, for instance, \cite[Ch.~5]{Kat1992})

Note that we may realize $A\coloneqq\big(\frac{a,b}{\Q}\big)$ as a $\Q$-subalgebra of $\mathrm{Mat}_2(\C)$ via  
\begin{equation}\label{Eq: generators}
		\Unit =
		\begin{pmatrix}
			1 & 0 \\ 0 & 1
		\end{pmatrix},\quad
		I = \sqrt a i
		\begin{pmatrix}
			0 & 1 \\ -1 & 0
		\end{pmatrix},\quad
		J=\sqrt b
		\begin{pmatrix}
			0 & 1 \\ 1 & 0
		\end{pmatrix},\quad
		K=\sqrt{ab}i 
		\begin{pmatrix}
			1 & 0 \\ 0 & -1
		\end{pmatrix}.
\end{equation}

The reduced norm of an element $Q = q_0\Unit + q_1I+ q_2J + q_3K\in A\subset \mathrm{Mat}_2(\C)$ is given by $\det(Q)$. Explicitly, we have 
\begin{equation*}
	Q = 
	\begin{pmatrix}
		q_0 + \sqrt{ab}iq_3 & \sqrt{a}iq_1 + \sqrt{b}q_2 \\ -\sqrt{a}iq_1 + \sqrt{b}q_2& q_0 - \sqrt{ab}iq_3 
	\end{pmatrix},
\end{equation*}
so 
\begin{equation*}
	\det(Q) = q_0^2 -aq_1^2-bq_2^2 + abq_3^2.
\end{equation*}
Note that if  $a$ and $b$ are positive and $\det(Q) = 1$, then $Q\in\SU(1,1)$. 

The \emph{standard order} in $A$ is given by 
\begin{equation*}
	\mc O = \{ Q = q_0 \Unit + q_1I + q_2J + q_3K : q_i\in \Z\}
\end{equation*}
We can now define the Fuchsian group associated with the standard order in $A$. It is given by the elements of $\mc O$ of unit reduced norm: 
\begin{equation*}
	\Gamma(A,\mc O)\coloneqq \{Q\in\mc{O}: \det(Q)  = 1\}\subset \SU(1,1).
\end{equation*}
It is known that $\Gamma(A,\mc O)$ constructed as above is a Fuchsian group. Moreover, the fact that $A$ is a division algebra guarantees that the quotient space $\CH^1/\Gamma(A,\mc O)$ is compact \cite[Ch.~5]{Kat1992}. In particular, $\Gamma(A,\mc O)$ is a co-compact Fuchsian group.

We shall see next that $\Gamma(A,\mc O)$ preserves a lattice in $\Heis_{5}$, as desired. In what follows, we shall use the following realization of $\Heis_{2n+1}$ for any $n\in\mathbb N$. Consider $\C^n$ with the Hermitian structure $(h,\bar\cdot)$ of signature $(1,n-1)$ and associated symplectic form $\omega= \mathrm{Im}(h)$. Then $\Heis_{2n+1}$ is the set $\C^n\times\R$ with multiplication law 
\begin{equation}\label{eq:Heis2n+1}
	(v,t)\cdot(v',t') = (v+v',t+t'+\tfrac{1}{2}\omega(v,w)),
\end{equation}
for any $v,v'\in\C^n, t,t'\in\R$
\begin{prop}\label{prop:HeisLattice}
	Let $b$ be a prime number and $a\in \N$ a quadratic non-residue modulo $b$. Equip $A=\big(\frac{a,b}{\Q}\big)$ with the standard order $\mc O$. Then there exists a $\Gamma(A,\mc O)$-invariant and co-compact lattice $\Gamma_2$ in the five-dimensional Heisenberg group. Thus, $\bar\Gamma_{a,b}\coloneqq \Gamma(A,\mc O)\ltimes\Gamma_2$ is a co-compact lattice in $\SU(1,1)\ltimes\Heis_5$.
\end{prop}
\begin{myproof}
	We write $e_1 = 
	\begin{psmallmatrix}
	1 \\ 0
	\end{psmallmatrix}\in \C^2$. Looking at \eqref{Eq: generators}, we see that 
	\begin{equation*}
		\mc O\cdot e_1 \coloneqq \Span_\Z\{e_1, Ie_1,Je_1,Ke_1\} 
		= \Span_{\Z}\{e_1, -\sqrt{a}ie_2,\sqrt{b}e_2,\sqrt{ab}ie_1\}
	\end{equation*}
	is a lattice in the real vector space $\C^2\cong \R^4$. Now we consider the subgroup of $\Heis_5$ generated by the lattice $\mc O\cdot e_1$. To see that the result is a lattice in the Heisenberg group, we need only check that the possible entries in the additional $\R$-factor are discrete. These are determined by the symplectic form corresponding to the Hermitian metric of signature $(1,1)$ on $\C^2$. By \eqref{eq:Heis2n+1}, the values that occur in the subgroup generated by a lattice $\Lambda\subset \C^2$ are half-integer multiples of the evaluation of this symplectic form on its generators. Thus, it suffices to determine these:
	\begin{gather*}
		\omega(Ie_1,e_1)=\omega(Je_1,e_1)=0=\omega(Ie_1,Ke_1)=\omega(Je_1,Ke_1)\\
		\omega(e_1,Ke_1)=\sqrt{ab}=\omega(I e_1,Je_1)
	\end{gather*} 
	Clearly, we only obtain multiples of $\sqrt{ab}$ by evaluation. This shows that the subgroup generated by this lattice in $\C^2$   is a co-compact lattice in $\Heis_5$.
\end{myproof}

We have now constructed infinitely many co-compact lattices $\bar\Gamma_{a,b}\subset\SU(1,1)\ltimes\Heis_5$, labeled by pairs $(a,b)\in \N^2$, where $b$ is prime and $a$ a quadratic non-residue modulo $b$. 

\begin{thm}\label{thm:mainn=2}
	Let $a,b\in\N$ be such that $b$ is prime and $a$ is a quadratic non-residue modulo $b$. Choose $c\geq 0$ such that $\tfrac{1}{2}\sqrt{ab}$ and $4\pi c$ are linearly dependent over $\Q$. Then $\bar\Gamma_{a,b}$ admits a co-compact sub-lattice $\tilde\Gamma_{a,b}$ such that $(\hat N_2^c,\hat g^c)/\tilde\Gamma_{a,b}$ is a complete quaternionic K\"ahler manifold diffeomorphic to $\R_{>0}\times K$, where $K$ is a compact and locally homogeneous seven-dimensional manifold.
\end{thm}
\begin{myproof}
	By \Cref{prop:G'structure} and \Cref{thm:SUltimesHaction} and since $\tfrac{1}{2}\sqrt{ab}$ and $4\pi c$ are linearly dependent over $\Q$, we see that $\bar\Gamma_{a,b}\cap Z(\Heis_5)+\mathcal F'$ is a discrete subgroup of $Z(\Heis_5) \cong \R$. It follows that $\bar\Gamma_{a,b}$ descends to a lattice $\hat\Gamma_{a,b}$ in $\SU(1,n-1)\ltimes (\Heis_{5}/\mathcal F')$ which acts effectively on $\hat N_2^c$. If $\hat\Gamma_{a,b}$ acts freely on $(\hat N_2^c,\hat g^c)$, we may directly take the quotient and obtain a quaternionic K\"ahler manifold with the required properties. We may think of $\hat N_2^c\cong\R_{>0}\times \bar N_{\rho_0}/\mathcal F'$, where $\bar N_{\rho_0} =\{\rho\equiv\rho_0>0\}\subset \bar N_2$, as a fiber bundle over $\R_{>0}$ with homogeneous fibers. The $\rho$-coordinate is preserved by all isometries, so $\hat N_2^c/\hat\Gamma_{a,b}$ is a fiber bundle over $\R_{>0}$ with locally homogeneous fibers which are moreover compact since $\hat\Gamma_{a,b}$ is co-compact.
	
	However, the constructed lattices $\hat\Gamma_{a,b}$ do not necessarily act freely on $\hat N_2^c$. Nevertheless their intersection with the stabilizer of a point, which is compact by \Cref{cor:compactstab}, is finite. Since a finite-index subgroup of a co-compact lattice is once again a co-compact lattice, it now suffices to find a finite-index subgroup of $\hat\Gamma_{a,b}$ which only intersects this finite group in the identity.
	
	The existence of such a subgroup is guaranteed by Selberg's lemma, which asserts that every finitely generated subgroup of $\mathrm{GL}(n,\C)$ admits a finite-index normal subgroup which is torsion-free. To apply this result, it remains to check that $\bar\Gamma_{a,b}$ is finitely generated. Observe that $\bar\Gamma_{a,b} = \Gamma(A,\mc O)\ltimes \Gamma_2$, where $\Gamma_2\subset\Heis_5$ is the lattice generated by $\{e_1,Ie_1,Je_1,Ke_1\}$, hence finitely generated. It therefore suffices to check that for the Fuchsian group $\Gamma(A,\mc O)$ is finitely generated. But $\Gamma(A,\mc O)$ is a lattice in a semi-simple Lie group and all such lattices are even finitely presented \cite[Ch.~4]{Mor2015}, so Selberg's lemma applies. Thus, we obtain a co-compact lattice $\tilde\Gamma_{a,b}$ which acts freely and isometrically on $(\hat N_2^c,\hat g^c)$ so that the corresponding quotient space possesses all the claimed properties. 
\end{myproof}

We now move on to describe a class of lattices $\Gamma\subset\SU(1,n-1)\ltimes\Heis_{2n+1}$ for arbitrary $n\geq 2$. These are of the form $\Gamma = \bar\Gamma_1\ltimes\Gamma_2$, where $\Gamma_2$ is a lattice in $\Heis_{2n+1}$ with the property that its normalizer $\bar \Gamma_1\subset\SU(1,n-1)$ is again a lattice. We consider lattices $\Gamma_2\subset\Heis_{2n+1}$ that are generated by lattices $\Lambda\subset\C^n$ satisfying a compatibility condition with the Hermitian structure $(h,\overline{\cdot})$ on $\C^n$ of signature $(1,n-1)$.
 
To state this condition, let $F$ be a totally imaginary quadratic number field, i.e.\ $F = \Q[i\sqrt{d}]$, where $d\in\N$ is a square-free integer. If we denote by $\mc O_F$ the ring of integers of $F$, then we have $\mc O_F = \Z[i\sqrt{d}]$, if $d \equiv 1,2\mod 4$ and $\mc O_F = \Z[\tfrac{1+i\sqrt{d}}{2}]$ if $d \equiv 3\mod 4$.

We call a lattice $\Lambda\subset\C^n$ that admits an action of the ring $\mc O_F$ of integers in $F$ \emph{compatible} \label{compatible:def_in_text} with $(h,\bar{\cdot})$ if $\bar\Lambda\subset\Lambda$, i.e.~$\Lambda$ is invariant under complex conjugation, and $h|_{\Lambda\times\Lambda}$ induces a sesquilinear form over $\mc O_F$. 

We remark that in the co-compact lattices $\bar{\Gamma}_{a,b} = \Gamma(A,\mathcal O)\ltimes \Gamma_2$ constructed above for $n=2$, the lattice $\Gamma_2\subset \Heis_5$ is generated by a lattice $\Lambda\subset \C^2$ satisfying this compatibility condition with $d = ab$, provided $d$ is square-free and $d\equiv 1,2\mod 4$. 

\begin{prop}
	Let $n\in \N$ be arbitrary and denote the generator of $Z(\Heis_{2n+1})\cong \R$ by $T$. Then, for any square-free $d\in \N$ with $d\equiv 1,2\mod 4$ there exists a lattice $\Gamma_2\subset \Heis_{2n+1}$ generated by a lattice $\Lambda\subset\C^n$ compatible with $(h,\bar\cdot)$ such that $\Gamma_2\cap Z(\Heis_{2n+1}) = \tfrac{1}{2}\sqrt{d}\Z T$ and the stabilizer $\bar\Gamma_1$ in $\SU(1,n-1)$ is also a lattice.  In particular, $\bar\Gamma\coloneqq \bar\Gamma_1\ltimes\Gamma_2\subset \SU(1,n-1)\ltimes\Heis_{2n+1}$ is a lattice. 
		
	If $n\geq 3$ then a lattice $\bar{\Gamma}\subset \SU(1,n-1)\ltimes\Heis_{2n+1}$ constructed in this way is never co-compact.
\end{prop}
\begin{myproof}
	The assertion is clear for $n=1$, so assume that $n\geq 2$. Let $d\in\mathbb{N}$ be square-free with $d\equiv 1,2\mod 4$, so that $F=\Q[i\sqrt{d}]$ is a totally imaginary quadratic number field with ring of integers $\mc O_F = \Z[i\sqrt{d}]$. 
	
	We start by constructing a suitable lattice $\Lambda_d\subset\C^n$, which is compatible with $(h,\bar{\cdot})$. Let $\{e_j: j=1,\dots,n\}$ be the standard basis of $\C^n$, which is orthonormal with respect to the Hermitian metric $h$ of signature $(1,n-1)$ on $\C^n$. Setting $f_j=i e_j$, we may consider the lattice $\Lambda_d\subset \R^{2n}\cong \C^n$ generated by $\{e_j,\sqrt{d} f_j: j=1,\dots,n\}$. Since the symplectic form associated to $h$ takes values in $\Z \cdot \sqrt{d}$ when evaluated on this lattice in $\C^n$, the subgroup $L_d$ of $\Heis_{2n+1}$ that this lattice generates is in fact a co-compact lattice. By construction, the above lattice $\Lambda_d$ is preserved by the natural action of $\mc O_F$. Moreover, $\Lambda_d$ is invariant under complex conjugation and $h|_{\Lambda\times\Lambda}$ induces a sesquilinear form over $\mc O_F$, so $\Lambda$ is compatible with $(h,\bar\cdot)$ in the sense defined above. 
	
	In this fashion, any square-free $d\in \N$ with $d\equiv 1,2\mod 4$ gives rise to a (co-compact) lattice $L_d\subset \Heis_{2n+1}$ which is generated by a lattice $\Lambda_d\subset \C^n$ compatible with $(h,\bar\cdot)$. 
	
	Recall that we view $\Heis_{2n+1}$ as the set $\C^n\times\R$ with multiplication law given by 
	$$(v,t)\cdot(v',t') = (v+v',t+t'+\tfrac{1}{2}\omega(v,w)),$$ where $\omega = \mathrm{Im}(h)$. Note that the Lie algebra $\heis_{2n+1}$ is the (real) vector space $\C^n\oplus\R$ with Lie bracket 
	\begin{equation*}
		[(v,t),(v',t')] = \omega(v,v').
	\end{equation*}
	In this setting the exponential map $\exp:\heis_{2n+1}\to\Heis_{2n+1}$ is just the identity on $\C^n\times\R$. 

	In general, we have the following description of the lattice $\Gamma_2\subset \Heis_{2n+1}$ generated by a lattice $\Lambda\subset\C^n$ compatible with $(h,\bar{\cdot})$. Let $rT$ be a generator for the image $\omega(\Lambda\times\Lambda)\subset\R$. Then we have 
	\begin{equation*}
		\Gamma_2 = \{(v,t)\in\C^n\times\R: v\in\Lambda, t\in\tfrac{1}{2}rT\Z\}.
	\end{equation*}
	This follows from the observation that if $v,v'\in\Lambda$ satisfy $\omega(v,v') = rT$, then we have in $\Heis_{2n+1}$: 
	\begin{equation*}
		(v,0)\cdot(v',0)\cdot(-(v+v'),0) = (0,\tfrac{1}{2}rT).
	\end{equation*}
	In particular, for the lattice $L_d$ constructed above, we find $L_d\cap Z(\Heis_{2n+1}) = \tfrac{1}{2}\sqrt{d}\Z T$.

	The set $\log\Gamma_2\subset\heis_{2n+1}= \C^n\oplus\R$ is closed under addition and equals the lattice $\Lambda\oplus \frac{1}{2}rT\Z$. Thus, the lattice $\Gamma_2$ generated by a compatible lattice $\Lambda\subset\C^n$ is a log-lattice. The lattice $\log \Gamma_2\subset\heis_{2n+1}$ induces the $\Q$-structure $(\heis_{2n+1})_{\Q} = \mathrm{span}_{\Q}(\log\Gamma_2)$ and the compatibility of $\Lambda$ ensures that $\SU(1,n-1)$ is defined over $\Q$ with respect to $(\heis_{2n+1})_{\Q}$. It then follows from \cite[Thm.~2.2]{MM1994} that the stabilizer $\bar\Gamma_1\subset\SU(1,n-1)$ of $\Gamma_2$ is a lattice in $\SU(1,n-1)$.  	We now take $\Gamma_2 = L_d$ and the above discussion shows that $\bar\Gamma_1 = \mathrm{Stab}_{\SU(1,n-1)}(\Gamma_2)$ is a lattice. We may then form the semi-direct product to obtain a lattice $\bar\Gamma =  \bar\Gamma_1 \ltimes \Gamma_2\subset \SU(1,n-1)\ltimes\Heis_{2n+1}$.

	If $n\geq 3$, and $\Lambda\subset \C^n$ is a lattice compatible with $(h,\bar\cdot)$ generating a lattice $\Gamma_2\subset\Heis_{2n+1}$, then the lattice $\bar\Gamma_1 = \mathrm{Stab}_{\SU(1,n-1)}(\Gamma_2)\subset \SU(1,n-1)$ contains a unipotent element and is therefore not co-compact by Godement's compactness criterion \cite[Prop.~5.3.1]{Mor2015}. To find such a unipotent element, we consider the indefinite rational quadratic form given by the real part of $h$ on the $2n$-dimensional $\Q$-vector space $\mathrm{span}_{\Q}(\Lambda)$. Since $n\geq 3$, i.e. $2n\geq 6$, we can apply Meyer's theorem \cite[Cor.~2 on p.43]{Ser73}  to find a vector $v \in\Lambda$ such that $h(v,v) = 0$. Furthermore, we can choose $w\in\Lambda$ such that $v$ and $w$ are linearly independent and $h(v,w) = 0$. Then $A\in\End(\C^n)$ given by $A = h(\cdot, v)w-h(\cdot,w)v$ is not zero, skew-hermitian with respect to $h$ and nilpotent, in particular $A\in \mf{su}(1,n-1)$. Thus, the unipotent element $\exp(A)\in\SU(1,n-1)$ stabilizes $\Lambda$, and hence $\Gamma_2$, i.e. $\exp(A)\in\Gamma_1$. For the lattices $L_d$ constructed above, we can take $v = e_1+e_2$ and $w = \sqrt{d}(f_1+f_2)$. This shows that in fact also for $n=2$ the stabilizer of $L_d$ in $\SU(1,n-1)$ is not co-compact.
\end{myproof}

Proceeding as in the proof of \Cref{thm:mainn=2}, we now easily obtain:

\begin{thm}\label{thm:finitevol}
	Let $n\in \N$ be arbitrary and let $d\in \N$ be square-free such that $d\equiv 1,2\mod 4$. Then there exists a lattice $\bar\Gamma\subset \SU(1,n-1)\ltimes\Heis_{2n+1}$ such that $\bar\Gamma\cap Z(\Heis_{2n+1})=\Z \cdot \frac{1}{2}\sqrt d T$, where $T$ is the generator of $Z(\Heis_{2n+1})\cong \R$. Choose $c\geq 0$ such that $4\pi c$ and $\tfrac{1}{2}\sqrt d$ are linearly dependent over $\Q$. Then $\bar\Gamma$ contains a lattice $\tilde\Gamma$ such that $(\hat N_n^c,g^c)/\tilde\Gamma$ is a complete quaternionic K\"ahler manifold diffeomorphic to $\R \times K$, where the fibers $\{t\} \times K$ are locally homogeneous and of finite volume.\proofclear
\end{thm}

\subsection{The volume of fiberwise quotients by lattices}

Given a lattice in $\SU(1,n-1)\ltimes\Heis_{2n+1}$, we have outlined above how to obtain a quaternionic K\"ahler manifold $\bar N_n/\Gamma$ which can be viewed as a fiber bundle with locally homogeneous fibers of finite volume. In the previous section we have shown how to obtain infinitely many lattices for any $n\in \N$, and constructed infinitely many co-compact examples in the case $n\leq 2$. We are now interested in the dependence of the volume of the fibers $\bar N_\rho/\Gamma\subset \bar N_n/\Gamma$ on $\rho$.

Recall the explicit expression for the one-loop deformed $c$-map metrics on $\bar N_n$ given in \eqref{eq:Bergmanmetric} and \eqref{eq:fullmetric}:
\begin{align*}
	g^c&=\frac{\rho+c}{\rho}g_{\CH^{n-1}}+\frac{1}{4\rho^2}\frac{\rho+2c}{\rho+c}\d \rho^2\\
	&\quad +\frac{1}{4\rho^2}\frac{\rho+c}{\rho+2c}
	\bigg(\d \tilde\phi 
	-4\Im\bigg(\bar w^0 \d w^0-\sum_{a=1}^{n-1}\bar w^a \d w^a\bigg)
	+\frac{2c}{1-\abs{X}^2}\Im \bigg(\sum_{a=1}^{n-1} \bar X^a \d X^a\bigg)\bigg)^2\\ 
	&\quad -\frac{2}{\rho}\bigg(\d w^0 \d \bar w^0-\sum_{a=1}^{n-1}\d w^a \d \bar w^a\bigg)
	+\frac{\rho+c}{\rho^2}\frac{4}{1-\abs{X}^2}	\Big|\d w^0+\sum_{a=1}^{n-1}X^a\d w^a\Big|^2,
\end{align*}
where
\[
	g_{\CH^{n-1}}=\frac{1}{1-\abs{X}^2}\Bigg(\sum_{a=1}^{n-1} \abs{\d X^a}^2+\frac{1}{1-\abs{X}^2}\bigg|\sum_{a=1}^{n-1}\bar X^a \d X^a\bigg|^2\Bigg).
\]
Ordering our coordinates as $(\rho,X,w,\tilde\phi)$, where we wrote $X = (X^1,\dots,X^{n-1})$ and $w = (w^0,\dots,w^{n-1})$,  the Gram matrix $G^c$ of $g^c$ takes on the following block form: 
\begin{equation*}
	G^c = \begin{pmatrix}
		g_{\rho\rho} & 0 & 0 & 0 \\ 0 & g_{XX} & g_{Xw} & g_{X\tilde\phi}\\ 0 & g_{Xw}^\top & g_{ww} & g_{w\tilde\phi}\\ 0 & g_{X\tilde\phi}^\top  & g_{w\tilde\phi}^\top  & g_{\tilde\phi\tilde\phi} 
	\end{pmatrix}
\end{equation*}
where we have coefficients $g_{\rho\rho},g_{\tilde\phi\tilde\phi}\in C^\infty(\bar N_n,\R)$, $g_{XX}\in C^\infty(\bar N_n,\mathrm{Mat}_{2n-2,2n-2}(\R))$, $g_{ww}\in C^\infty(\bar N_n,\mathrm{Mat}_{2n,2n}(\R))$, $g_{Xw}\in C^\infty(\bar N_n,\mathrm{Mat}_{2n-2,2n}(\R))$, $g_{X\tilde\phi}\in C^\infty(\bar N_n,\R^{2n-2})$ and $g_{w\tilde\phi}\in C^\infty(\bar N_n,\R^{2n})$.

Since the level sets of $\rho$ are locally homogeneous of finite volume, it 
suffices to compute the Gram matrix and its determinant at a point $(\rho ,p_0)$ where $X=0$ and $w=0$. At such a point the metric $g^c$ evaluates to
\begin{align*}
	g^c(\rho,p_0) &=\frac{\rho+c}{\rho}\Bigg(\sum_{a=1}^{n-1} \abs{\d X^a}^2\Bigg)+\frac{1}{4\rho^2}\frac{\rho+2c}{\rho+c}\d \rho^2 +\frac{1}{4\rho^2}\frac{\rho+c}{\rho+2c}
	\d \tilde\phi^2\\ 
	&\quad+\frac{2}{\rho}\bigg(\frac{\rho+2c}{\rho}\abs{\d w^0}^2 +\sum_{a=1}^{n-1}\abs{\d w^a}^2\bigg).
\end{align*} 
In particular, the off-diagonal components $g_{Xw}, g_{X\tilde\phi}, g_{w\tilde\phi}$ of $G^c$ all vanish at  $(\rho,p_0)$ and we have
\[
	g_{\rho\rho}(\rho,p_0) = \frac{1}{4\rho^2}\frac{\rho+2c}{\rho +c},\qquad 
	g_{XX}(\rho,p_0) =  \frac{\rho + c}{\rho}\Unit_{2n-2}
\]
and 
\[
g_{ww}(\rho,p_0) = \frac{2}{\rho} \begin{pmatrix}
\frac{\rho + 2c}{\rho}\Unit_2 & 0 \\ 0 & \Unit_{2n-2}
\end{pmatrix} \qquad	
	g_{\tilde\phi\tilde\phi}(\rho,p_0) = \frac{1}{4\rho^2}\frac{\rho+c}{\rho+2c}.
\]
It follows that 
\[
\det G^c(\rho,p_0) = \frac{2^{2n-4}}{\rho^{2n+4}}\left(\frac{\rho+c}{\rho}\right)^{2n-2}\left(\frac{\rho+2c}{\rho}\right)^2. 
\]
Thus, $f=\sqrt{\det G^c} = \frac{1}{\rho^{n+2}}\left(\frac{\rho+c}{\rho}\right)^{n-1}\left(\frac{\rho+2c}{\rho}\right) f_{\mathrm{inv}}$,  where $f_{\mathrm{inv}}$
is the volume density of the unique invariant volume form $f_{\mathrm{inv}}\d X\,\d w\,\d \tilde\phi$ in the fiber normalized by  $f_{\mathrm{inv}}(p_0)=2^{n-2}$. Here $\d X\,\d w\,\d \tilde\phi$ stands for the Lebesgue measure associated with the coordinate system.

Note that we have for any $c\geq 0$ and $\rho >0$ the inequality
\begin{equation}\label{eq:ineq_rho0}
f=  \frac{1}{\rho^{n+2}}\left(1+\frac{c}{\rho}\right)^{n-1}\left(1+\frac{2c}{\rho}\right) f_{\mathrm{inv}} \geq  \frac{f_{\mathrm{inv}}}{\rho^{n+2}}.  
\end{equation}
On the other hand, if  $\rho\geq \rho_0>0$ then we have  for any $c\geq 0$, 
\begin{equation}\label{eq:ineq_rhoinfty}
f = \frac{1}{\rho^{n+2}}\left(1+\frac{c}{\rho}\right)^{n-1}\left(1+\frac{2c}{\rho}\right) f_{\mathrm{inv}} \leq  C(\rho_0)\frac{f_{\mathrm{inv}}}{\rho^{n+2}}, 
\end{equation}
where $C(\rho_0)  = \left(1+\frac{c}{\rho_0}\right)^{n-1}\left(1+\frac{2c}{\rho_0}\right)>0$. 

Fix $\rho_0>0$ and let $D$ denote a fundamental domain for the action of $\Gamma$ on the fiber $\bar N_{\rho_0}$. Then we deduce from \eqref{eq:ineq_rhoinfty} that there exists some $C_1>0$ such that
\begin{align*}
	\vol_{\bar N_{\rho\geq \rho_0}/\Gamma}
	&\coloneqq \lim_{\rho\to \infty}\vol_{\bar N_{[\rho_0,\rho]}/\Gamma}
	=\lim_{\rho\to\infty}\int_{D\times[\rho_0,\rho]}f \d X\,\d w\,\d \tilde\phi\,\d \rho\\
	&\leq C_1 \int_{\rho_0}^\infty \frac{\d \rho}{\rho^{n+2}
	=\frac{C_1}{(n+1)\rho_0^{n+1}}}
\end{align*}
Similarly, we obtain from \eqref{eq:ineq_rho0} the existence of a constant $C_2>0$ such that
\begin{align*}
 	\vol_{\bar N_{[\rho_1,\rho_0]}/\Gamma}&=  \int_{D\times[\rho_1,\rho_0]}f \d X\,\d w\,\d \tilde\phi\,\d \rho\\
 	& \geq C_2 \int_{\rho_1}^{\rho_0}\frac{\d \rho}{\rho^{n+2}}
 	=\frac{C_2}{(n+1)}\bigg(\frac{1}{\rho^{n+1}_1}-\frac{1}{\rho_0^{n+1}}\bigg)
\end{align*}
for any $0<\rho_1 <\rho_0$. We summarize this discussion as follows.

\begin{thm}\label{thm:volume}
	Let $c\geq 0$ and consider the quaternionic K\"ahler manifold $(\bar N_n/\Gamma, g^c)$. Then for each $\rho_0\in \R_{>0}$ we have 
	\begin{equation*}
		\vol(\bar N_{\rho\geq\rho_0}/\Gamma)<\infty,
	\end{equation*}
	and
	\begin{equation*}
		\vol(\bar N_{[\rho,\rho_0]}/\Gamma)\to\infty,\quad \text{as $\rho\to 0$}.\vspace{-0.7cm}
	\end{equation*}
	\proofclear
\end{thm}
\begin{rem}
We may write $\frac{f}{f_{\mathrm{inv}}}  = \frac{1}{\rho^{n+2}}\left(\frac{\rho+c}{\rho}\right)^{n-1}\left(\frac{\rho+2c}{\rho}\right) = \frac{1}{\rho^{n+2}}P\left(\frac{c}{\rho}\right)$, where $P$ is a polynomial of degree $n$ with positive coefficients and constant term equal to $1$. From this one may obtain the explicit formula
\begin{equation*}
	\vol_{\bar N_{[\rho_1,\rho_0]}/\Gamma} =  V(D) \int_{\rho_1}^{\rho_0}\frac{1}{\rho^{n+2}}P\left(\frac{c}{\rho}\right)\d\rho, \qquad V(D) =\int_D f_{\mathrm{inv}} \d X\,\d w\,\d \tilde\phi.
\end{equation*} 
Analyzing this expression gives the leading asymptotics 
\[ \vol_{\bar N_{\rho\geq \rho_0}/\Gamma} = k\rho_0^{-n-1} + O(\rho_0^{-n-2}),\qquad k=  \frac{V(D)}{n+1}  \]
when $\rho_0\rightarrow \infty$, while the asymptotics of  $\vol_{\bar N_{[\rho_1,\rho_0]}/\Gamma}$ for $\rho_1 \rightarrow 0$ depends on $c$. 
For $c=0$ the volume is simply 
\[ \vol_{\bar N_{[\rho_1,\rho_0]}/\Gamma} = k(\rho_1^{-n-1}-\rho_0^{-n-1}) \sim k\rho_1^{-n-1},\quad \rho_1\rightarrow  0,  \]
whereas for $c>0$ it grows like $k_1\rho_1^{-2n-1}$:
\[ \rho_1^{2n+1}\vol_{\bar N_{[\rho_1,\rho_0]}/\Gamma} = k_1 + O(\rho_1),\quad k_1=\frac{2c^n}{2n+1}k    \]

	In the case $n=1$ there are no $X$-coordinates and the above computations simplify considerably. The Gram matrix $G^c$ takes values in $\mathrm{Mat}_{4,4}(\R)$ and $\frac{i}{2}\d w\wedge \d\bar w\wedge \d\tilde\phi$ is a left-invariant volume form on $\Heis_{3}$. The determinant $\det G^c$ is therefore a function of $\rho$ alone, i.e.\ $f_{\mathrm{inv}}$ is constant. 
\end{rem}

\subsection*{Conflict of interest statement}
Vicente Cort\'es is managing editor of ``Abhandlungen aus dem Mathematischen Seminar der Universit\"at Hamburg'' as well as editor of ``Beitr\"age in Algebra und Geometrie'' and  ``Differential Geometry and its Applications.'' Funding information for this paper is included under ``Acknowledgements'' after the introduction.  

On behalf of all authors the corresponding author states that here is no conflict of interest.

\subsection*{Data availability statement}
There are no data for this article.

\end{document}